\documentclass{siamart220329}
\usepackage{amsmath, amssymb}
\usepackage{mathtools}
\usepackage{hyperref}
\usepackage{verbatim}
\usepackage{tikz}
\usepackage[mathscr]{euscript}
\usepackage{cite}
\usepackage[normalem]{ulem}
\usetikzlibrary{cd}

\providecommand{\Div}{\operatorname{div}}          

\providecommand{\Dim}{\operatorname{dim}}            
\providecommand{\dim}{\Dim}

\providecommand{\argmin}{\operatorname*{argmin}}  

\providecommand{\Ce}{{\cal E}}

\providecommand{\Ch}{{\cal H}}

\providecommand{\Ck}{{\cal K}}

\providecommand{\bbN}{\mathbb{N}}

\providecommand{\bbR}{\mathbb{R}}

\newcommand{\ds}{\,\mathrm{d}s}

\newcommand{\dx}{\,\mathrm{d}x}

\newcommand*{\jump}[1]{\left[\!\left[{#1}\right]\!\right]\!}
\newcommand*{\avg}[1]{\left\{{#1}\right\}\!}
\newcommand*{\jumpe}[1]{\left[\!\left[{#1}\right]\!\right]\!|_{e}}
\newcommand*{\avge}[1]{\left\{{#1}\right\}\!|_{e}}
\newcommand*{\upwind}[1]{\left\lfloor {#1} \right\rfloor}

\newsiamremark{remark}{Remark}
\newsiamremark{example}{Example}
\newsiamremark{numexp}{Numerical Experiment}

\title{Level-set shape optimization \\ via polytopic discontinuous Galerkin
methods \thanks{\funding{The work of R.~S.~Fernandes is supported by the EPSRC
DTP grant EP/V520172/1. E.~H.~Georgoulis wishes to acknowledge the financial
support of EPSRC (grant number EP/W005840/2) and of The Leverhulme Trust (grant
number RPG-2021-238).}}} \author{Raphael S. Fernandes\thanks{School of
Computing and Mathematical Sciences, University of Leicester, Leicester LE1
7RH, United Kingdom (\email{rsf14@leicester.ac.uk})} \and Emmanuil
H.~Georgoulis\thanks{ The Maxwell Institute for Mathematical Sciences \&
Department of Mathematics, School of Mathematical and Computer Sciences,
Heriot-Watt University, Edinburgh EH14 4AS, United Kingdom\and Department of
Mathematics, School of Applied Mathematics and Physical Sciences, National
Technical University of Athens, Zografou Campus, Athens 15780, Greece\and
IACM-FORTH, Crete 70013, Greece (\email{e.georgoulis@hw.ac.uk})} \and Alberto
Paganini\thanks{School of Computing and Mathematical Sciences, University of
Leicester, Leicester LE1 7RH, United Kingdom
(\email{a.paganini@leicester.ac.uk})}}

\begin{document}
\maketitle

\begin{abstract}
We introduce a new level-set shape optimization approach based on polytopic
(i.e., polygonal in two and polyhedral in three spatial dimensions)
discontinuous Galerkin methods. The approach benefits from the geometric mesh
flexibility	of polytopic discontinuous Galerkin methods to resolve the
zero-level set accurately and efficiently. Additionally, we employ suitable
Runge-Kutta discontinuous Galerkin methods to update the level-set function on
a fine underlying simplicial mesh. We discuss the construction and
implementation of the approach, explaining how to modify shape derivate
formulas to compute consistent shape gradient approximations using
discontinuous Galerkin methods, and how to recover dG functions into smoother
ones. Numerical experiments on unconstrained and PDE-constrained test cases
evidence the good properties of the proposed methodology.
\end{abstract}

\section{Introduction}

The level-set method is a very popular shape optimization methodology based on
the implicit representation of domain boundaries (shapes) as zero-level sets of
a scalar function, which is often referred to simply as the \emph{level-set
function}. The level-set method has been successfully applied in a large
variety of scenarios including structural mechanics, architecture, and
multi-physics systems \cite{AlDaJo21, DaFaMiAlCoEs17, Fe19}, among others. A
notable feature of the level-set shape optimization methodology, which
certainly contributed to its widespread use, is that this method is capable of
performing also topological changes such as shape splitting or merging \cite{AlDaJo21}.
In particular, topological changes can occur even when level-set function
updates are driven only by shape derivatives, which encode information of
diffeomorphic (and thus non-topological) domain changes \cite{DeZo11}. Of
course, level-set optimization methods based on topological derivatives are
also available \cite{BlSt23, AmHe06}.

The level-set shape optimization methodology has been described in review and
introductory articles \cite{AlDaJo21, La18}, and it has been implemented in
open-source libraries \cite{WeMaMaBaCh24, AlGr24, letop}. Although many
different variants exist, the method typically involves solving a
\emph{shape-gradient equation} to extract a descent direction from the shape
derivative, and a \emph{level-set equation} to update the level-set function.
The shape-gradient equation is an elliptic boundary value problem with the shape
derivative as right-hand side. The level-set equation is a transport partial differential equation (PDE)
with non constant velocity (specifically, the shape gradient multiplied by -1).
If the shape optimization problem is constrained to partial differential
equations, then determining the shape derivative also requires solving
\emph{state} and \emph{adjoint equations} \cite{DeZo11}. In some instances, e.g.,
 in self-adjoint problems (but not only), the solution to the adjoint equation
can be obtained by rescaling the solution to the state equation \cite{AlDaJo21}.

Solving PDE-constrained shape optimization problems with the level-set method
generally requires three ingredients. The first one is a numerical method to
solve the level-set equation, which is typically stated on a bounded domain containing  
shapes,
often referred to as the \emph{hold-all domain}. The second ingredient is a numerical
method to solve the shape-gradient equation, whose right-hand side is supported
on the domain implicitly represented by the level-set function. The third one is
a numerical method to solve the state and the adjoint equations, which are
defined on the domain implicitly represented by the level-set function. 

The level-set equation is commonly solved using finite difference schemes on
Cartesian grids \cite{La18}; these schemes typically require level-set
reinitialization to ensure
stability \cite{MaReCh06}. Another popular approach is to add diffusion to the
level-set equation (in the spirit of vanishing viscosity \cite{Ev22}) and, 
subsequently, employ
standard conforming piecewise-linear finite elements \cite{letop, OtYaIzNi15}.
In \cite{BuDaFr12}, the authors introduce an alternative method based on the
method of characteristics and simplicial mesh adaptation, which is applied to
level-set shape optimization in \cite{AlDaFr14}.

Solving the shape gradient, state, and adjoint equations is more challenging
because these are partly or entirely defined on the domain implicitly
represented by the level-set function. Since this is usually not resolved by the
hold-all domain mesh, naively employing standard finite element methods leads to
substantial discretization errors. Numerous approaches have been developed to
address this issue. These can be classified as either being ``unfitted'' or
``fitted'' depending on how accurately they resolve the zero-level set.

Unfitted methods employ techniques to perform all computations on the hold-all
domain mesh. A very popular unfitted method is the so-called ``Ersatz material''
approach \cite{AlDaJo21}, which recasts the state equation onto the
hold-all domain by extending the equation's coefficients in a suitable fashion. In this
case, the shape gradient and the adjoint equations are also naturally
reformulated on the hold-all domain. When applicable, this approach allows
computing numerical solutions using standard finite element methods. However,
the main drawback of the Ersatz material approach is that the zero-level set
is not resolved accurately. Beside limiting computational accuracy, this
drawback disqualifies the Ersatz material approach for problems where a precise description of
the zero-level set is required as, for instance, in interface problems
\cite{AlDaFr14, Pa15}. These issues can be addressed by employing more
sophisticated unfitted method like the CutFEM
\cite{BuElHaLaLa19,BeWaBe18}, the XFEM \cite{DuVaJaFl06, ViMa14, MaMa14},
or the finite cell method \cite{PaDueRa12, MaThHiBa24}, among others.

Fitted methods employ techniques to accurately represent the implicitly defined
domain. Classical examples are based on mesh adaptation \cite{YaShTsKaSaNi11},
global re-meshing \cite{XiShShYu12}, or local re-meshing with mesh adaptation
\cite{AlDaFr14}. These approaches are effective, but can require careful
case-by-case treatment and mesh management to prevent the mesh quality from
deteriorating. In this work, we present a new fitted level-set shape
optimization method that drastically simplifies the mesh management step by
employing discontinuous Galerkin finite elements on polytopic meshes.

Polytopic meshes have been proved effective in shape and topology optimization 
thanks to their increased geometric flexibility. Polytopic Voronoi meshes with
low-order conforming finite elements \cite{TaPaPeMe10, TaPaPeMe12, HoMePe18} or
virtual elements \cite{AnBrScVe17, GaPaDuMe15, TrNgBuNg23} have been applied
successfully in combination with the ``Solid Isotropic Material with
Penalization'' (SIMP) model. In the SIMP model, the state equation is extended
to the hold-all domain similarly to the Ersatz material approach.  Polytopic
Voronoi meshes with low-order conforming finite elements have also been applied
to the level-set shape optimization method with ``Ersatz material'' approach
\cite{WePa20}. 

Polytopic meshes also greatly simplify fitted methods. In
\cite{NgKi19, NgSoKi22}, the authors discretize the level-set function using
standard finite elements on hexahedral meshes, truncate the hexahedra with the
zero-level set to generate a fitted polytopic mesh, employ modified trilinear
and moving-least squares basis functions to compute the state, and solve the
level-set equation using ``essentially non-oscillatory'' (ENO) schemes. In
\cite{FeMiPaPeVeAn24}, the authors present a level-set shape optimization based
on the topological derivative and generate fitted polytopic meshes by
intersecting an initial Voronoi tessellation with the zero-level set. The state
equation is solved using the symmetric Interior Penalty discontinuous Galerkin
(IPdG) method on the polytopic mesh. The level-set equation is regularized by
adding a small regularizing amount of diffusion and is solved on the polygonal
Voronoi mesh using a symmetric IPdG method. This approach does not include the
equivalent of the shape gradient equation because the topological derivative is
applied pointwise. {A different approach is presented in \cite{DaLeOu24},
where the authors combine
polytopic meshes induced by a Laguerre diagram with the virtual element method
and realized a shape and topology optimization method by acting directly on
the Laguerre diagram parameters.}

These works demonstrated some of the advantages polytopic meshes
offer in the context of shape optimization such as higher element anisotropy
or facilitated local refinement. However, polytopic
methods have much more to offer. The last ten years or so have seen the
development of a rich theory (supported by efficient implementations) of
discontinuous Galerkin (dG) methods on general, possibly curved, polytopic
elements \cite{CaDoGeHo14,CaDoGeHo16,CaDoGeHo17,DoGeKa21, CaDoGe21}.
In particular, the element-wise discontinuous nature of the approximations
spaces permitted by dG methods allows for extremely general element shapes.
{  The ability to use elements of essentially arbitrary shape with the same number of local numerical degrees of freedom as those used on simplicial elements \cite{CaDoGeHo14,CaDoGeHo16,CaDoGeHo17} effectively decouples the element shape from its approximation capabilities. This is extremely attractive in the context of shape optimization, since coarse polytopic elements can accurately represent complex shapes, in contrast to what is possible with coarse simplicial meshes. Moreover, hanging nodes are naturally embedded in the polytopic mesh paradigm.}

In this
work, we build on this enhanced geometric flexibility. Specifically, we propose a
\emph{polytopic discontinuous Galerkin framework} for fitted level-set shape
optimization based on agglomeration of elements of a fine underlying simplicial mesh.
The approach allows decoupling
the discretization of the level-set equation, the discretization of the shape
gradient equation, and the discretization of the state (and adjoint) equations.
This approach streamlines and facilitates the construction of coarse agglomerated
polytopic elements that \emph{fit the zero-level set to great accuracy} and that, in
combination with higher-order dG methods, allow the rapid solution of the shape
gradient and of the state equations. Transferring information between a coarse
polytopic mesh and the underlying fine simplicial mesh is also straightforward.
This enables solving the level-set equation on the fine underlying simplicial
mesh using standard, explicit, and stable Runge-Kutta dG (RKdG) methods, to a
computational cost comparable to that of solving the shape
gradient equation. 
We note that in \cite{MaReCh06} the authors evidence that using RKdG methods for
the level-set
equation in divergence form is very efficient, accurate, stable. Additionally,
in contrast to Hamilton-Jacobi ENO methods, RKdG methods do not require non-standard
procedures such as reinitialization. In this work, we apply RKdG methods directly
to the level-set equation in non-divergence form, without any observed
detrimental effect to the stability of the method.

The remainder of this work is organized as follows. In \cref{sec:shapeopttheory},
we briefly
review the level-set shape optimization method. In \cref{sec:dGlevelset}, we
discuss the discontinuous Galerkin approximation of the shape gradient, the
state, and the level-set equation, paying extra care to the consistent
evaluation of the shape derivative on discontinous test functions. In
\cref{sec:Numericalexperiments}, we present a selection of numerical experiments
that demonstrate the efficacy of the new method.  We consider both unconstrained
and PDE-constrained shape optimization test cases and use a steepest-descent
optimization algorithm with backtracking and Armijo stopping rule. The method is
implemented in Python and the code is available on the Zenodo archive at
\cite{FeGePa24}. Finally, in \cref{sec:conclusions} we draw conclusions and
outline future research directions.

\section{The level-set shape optimization method}\label{sec:shapeopttheory}
We review briefly the level-set shape
optimization method; for an extensive introduction to the topic; we refer to the educational articles \cite{AlDaJo21, La18}.

\begin{figure}[htb!]
\centering
\includegraphics[width=0.6\linewidth]{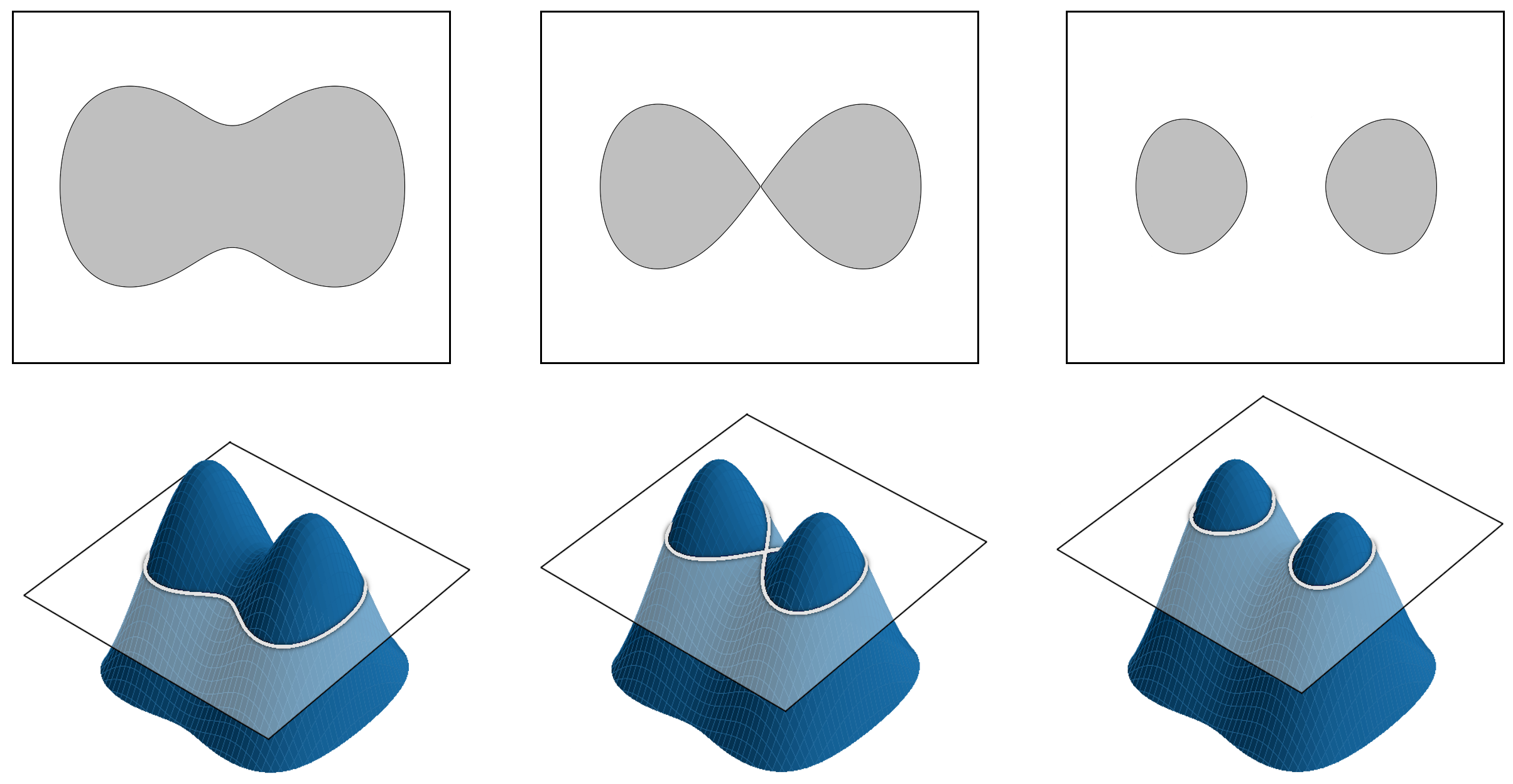}
\caption{Surf plot of a function $\phi$ and zero level sets of $\phi+c$, where
$c$ is a negative constant (decreasing from left to right). Note that the shape
and topology of $\Omega_{\phi+c}$ changes depending on the value of $c$.}
\label{fig:levelsetdomain}
\end{figure}

Let $D\subset\bbR^d$, $d\in\mathbb{N}$, be a fixed and Lipschitz hold-all
domain. The level-set shape optimization method is based on the idea that to
each continuous function $\phi:D\to \bbR$ we can associate a domain
\begin{equation*}
\Omega_\phi \coloneqq \{x\in D\,: \phi(x) < 0\}\,.
\end{equation*}
This identification allows modifying both the shape and the topology of a domain
$\Omega_\phi$ by acting on the defining function $\phi$. For example,
\cref{fig:levelsetdomain} shows that adding a negative constant to a
function $\phi$ can lead to changes in the shape and topology of $\Omega_\phi$.

Let $W^{1,\infty}(D,\bbR^d)$ denote the space of Lipschitz continuous functions
with values in $\bbR^d$. For a fixed function $V\in W^{1,\infty}(D,\bbR^d)$ with
$V=0$ on $\partial D$ and for $\epsilon>0$ sufficiently small, we consider the
family of bi-Lipschitz bijections $\{T_t:D\to D\}_{0\leq t\leq \epsilon}$
defined by $T_t(x) = x + tV(x)$ for all $x$ in $D$. By applying $T_t$ to
$\Omega_\phi$, we can create a family of domains $\Omega_\phi (t)\coloneqq
T_t(\Omega_\phi)$. In turn, we can create a family of functions
$\{\varphi(t):D\to\bbR\}_{0\leq t\leq \epsilon}$ such that $\Omega_{\varphi(t)}
= \Omega_{\phi}(t)$ by simply defining $\varphi(t, T_t(x)) = \phi(x)$ for all
$x$ in $D$ and all $0\leq t\leq \epsilon$. We henceforth refer to the variable $t$
as ``time''. By computing the total 
derivative of $\varphi(t,T_t)$ with respect
to $t$, we observe that $\varphi$ is characterized by the hyperbolic initial
value problem commonly known as
\emph{level-set equation}: 
\begin{equation} \label{eq:levelseteqn}
\partial_t \varphi + V\cdot \nabla\varphi = 0\quad \text{for }t\in(0,\epsilon)\,,
\quad \varphi(0) = \phi\,.
\end{equation}

Let $J:D\supset\Omega\mapsto J(\Omega)\in\bbR$ denote a sufficiently regular
objective function. Using the level-set characterization \cref{eq:levelseteqn},
we can formulate a level-set shape optimization method to tackle the shape
optimization problem:
\begin{equation*}
\text{Find }\Omega\subset D \text{ such that }\Omega\in\argmin_{A\subset D} J(A)\,.
\end{equation*}
More precisely, given an initial guess $\Omega_\phi$, we can generate a sequence
of functions $\phi^k$ such that $J(\Omega_{\phi^k})$ is decreasing by setting
$\phi^0=\phi$ and $\phi^{k+1} \coloneqq \varphi^k(t^k)$ for $k\geq 0$, where $\varphi^k$
is the solution to
\begin{equation}\label{eq:levelsetmethod}
\partial_t \varphi^k - \nabla J(\Omega_{\phi^k})\cdot \nabla\varphi^k = 0
\quad \text{for }t\in(0,\epsilon_k)\,, \quad \varphi^k(0) = \phi^{k}\,,
\end{equation}
$\nabla J(\Omega_{\phi^k})$ is the shape gradient of $J$ at $\Omega_{\phi^k}$
with respect to a suitable inner product, and $\epsilon_k>0$ and $t^k\in
(0,\epsilon_k)$ are suitable final and stopping times, respectively. Of course,
the optimization algorithm is stopped if $\nabla J(\Omega_{\phi^k})=0$, in which case
$\phi^k$ is a stationary point.
We recall
that the shape derivative $J$ at $\Omega$ in the direction $V\in
W^{1,\infty}(D,\bbR^d)$ is given by
\begin{equation*}
dJ(\Omega; V)\coloneqq \lim_{t\to 0} \frac{J(T_t(\Omega))-J(\Omega)}{t}\,.
\end{equation*}
Let $\big(\Ch, (\cdot,\cdot)_{\Ch}\big)$ be a Hilbert space of vector-valued
functions defined over $D$ such that $W^{1,\infty}(D,\bbR^d)$ is dense in $\Ch$.
In practice, a popular choice is to select $\Ch=H^1_0(D,\bbR^d)$, the space of 
$\bbR^d$-valued functions with components in $H^1_0(D)$, either endowed
with the standard $H^1$- or the $H^1_0$-inner products, or an inner product
stemming from the linear elasticity equations. Then, the shape gradient $\nabla
J(\Omega)\in \Ch$ of $J$ at $\Omega$ with respect to $\Ch$ is the solution to
\begin{equation}\label{eq:shapegradienteqn}
(\nabla J(\Omega), V)_\Ch = dJ(\Omega; V)\quad \text{for all } V\in \Ch\,.
\end{equation}

For existence and uniqueness of a weak solution to \cref{eq:levelseteqn}, it is
sufficient to assume that $V\in W^{1,\infty}(D)$ and, in particular, the same is
required for $V= -\nabla J(\Omega_{\phi^k})$. Since $\nabla J(\Omega_{\phi^k})$
will be approximated by a Galerkin method over piecewise polynomial spaces,
these approximations to $\nabla J(\Omega_{\phi^k})$ will be bounded in $D$ with
piecewise bounded derivatives. Therefore, for existence and uniqueness of
solution to \cref{eq:levelsetmethod}, it is sufficient to ensure the global
continuity of the approximations to $\nabla J(\Omega_{\phi^k})$  inserted in
\cref{eq:levelsetmethod}.

\begin{example}\label{ex:levelsetobjective}
As a driving example, let $f\in W^{1,1}(D)$ be a given scalar function, and let 
$J(\Omega)\coloneqq \int_\Omega f\dx$. Then, the shape derivative 
of $J$ at $\Omega\subset D$ in the direction $V\in W^{1,\infty}(D,\bbR^d)$ is given by
\begin{equation*}
dJ(\Omega; V) = \int_\Omega \big(\nabla f\cdot V + f\Div V\big)\dx\,.
\end{equation*}
\end{example}

\begin{remark}
As described in \cite{DeHeHi22}, it is also possible to compute steepest descent
directions in $W^{1,\infty}(D,\bbR^d)$. These descent directions can be used to drive the evolution
of the level-set function $\varphi^k$ in \cref{eq:levelsetmethod}
instead of using the negative shape gradient $-\nabla J(\Omega_{\phi^k})$. To the best of our understanding,
computing such descent directions may be computationally more involved than solving
\cref{eq:shapegradienteqn}. Although we are not aware
of their use in level-set based shape optimization, in other frameworks it has been
observed that $W^{1,\infty}(D,\bbR^d)$-steepest direction based methods can outperform Hilbert space-based
counterparts, especially when optimizing domains with corners. Since a comparison between different choices of
descent directions is not the focus of this work, for simplicity we compute descent directions by solving \cref{eq:shapegradienteqn}.
We highlight that, as presented below, the shape gradient equation
\cref{eq:shapegradienteqn} can be solved efficiently using
polytopic meshes due to their reduction of computational complexity.
\end{remark}

\section{Numerical level-set shape optimization using dG methods}
\label{sec:dGlevelset}

The numerical realization of the level-set shape optimization method
\cref{eq:levelsetmethod} requires, at the very least, the discretization of the
level set equation \cref{eq:levelseteqn} and of the shape gradient equation
\cref{eq:shapegradienteqn}. A very popular approach is to employ finite
differences to solve \cref{eq:levelseteqn} in combination with a standard
conforming Galerkin method to solve \cref{eq:shapegradienteqn}. Instead, here
we develop a \emph{polytopic discontinuous Galerkin method} framework.
The so-called \emph{polytopic} (i.e., polygonal for $d=2$ or polyhedral for
$d=3$) discontinuous Galerkin methods have been developed to be able to operate
on meshes comprising essentially arbitrary element shapes. In
\cite{CaDoGeHo14,CaDoGeHo16,CaDoGeHo17, CaDoGe21}, it has been shown
that it is possible to construct \emph{stable}
discontinuous Galerkin methods on such extreme mesh scenarios. A pertinent
feature of these methods is  that the cardinality of the local elemental degrees
of freedom are \emph{independent} of the element shape, that is element geometry
and number of faces/edges/vertices. This feature is naturally very attractive
for shape optimization frameworks, whereby it is extremely important to strike a
balance between total computational complexity, by using as coarse mesh as
possible, and geometric feature resolution. As we shall see below, by allowing
meshes with elements with essentially arbitrary shapes, we can follow zero-level
sets with extremely high fidelity without requiring highly locally refined
meshes in the zero-level set neighborhoods.

Before proceeding furher, we recall
some standard notation to introduce the discontinuous Galerkin method.
To this end, let $\mathscr{T}$ be a partition of the hold-all domain $D$ into
general polytopic elements. We consider the broken Sobolev space
\begin{equation*}
{H^1(\mathscr{T})\coloneqq
\{v\in L^1(D)\,: v\vert_T\in H^1(T) \text{ for all } T\in\mathscr{T}\},}
\end{equation*}
along with the broken gradient defined by $(\nabla_h v)|_T:=\nabla (v|_T)$,
$T\in\mathcal{T}$. For a generic element $T\in \mathscr{T}$, let ${P}_p(T)$
denote the space of tensor-product polynomials defined over $T$. We
consider the element-wise (discontinuous) polynomial space
\begin{equation}
S^p_{\mathscr{T}} := \{v \in L^2(D)\,:v|_T \in {P}_p(T)
\text{ for all } T\in\mathscr{T}\}\,.
\end{equation}
and we denote by $[S^p_{\mathscr{T}}]^d$ the space of vector-valued functions
whose $d$ components lie in $S^p_{\mathscr{T}}$. For later use we also consider the
restriction $S^p_{\mathscr{T},\Omega}$ of $S^p_{\mathscr{T}}$ onto the domain
$\Omega\subset D$ given by
\begin{equation}
S^p_{\mathscr{T},\Omega} := \{v \in L^2(\Omega)\,:v|_{T\cap\Omega} \in {P}_p(T\cap\Omega) \text{ for all } T\in\mathscr{T}\}\,,
\end{equation}
and we denote the standard, orthogonal $L^2(\Omega)$-projection
$\Pi_\Omega:L^2(\Omega)\to S^p_{\mathscr{T},\Omega}$ onto
$S^p_{\mathscr{T},\Omega}$,  defined for every $v\in L^2(\Omega)$ by
\begin{equation}\label{eq:L2proj}
\int_\Omega (v-\Pi_\Omega v) w \dx=0 \quad\text{for all } w\in S^p_{\mathscr{T},\Omega}. 
\end{equation}

Finally, we recall the standard definitions of jump and averages of a function
in $S^p_{\mathscr{T}}$. To this end, we introduce the notation $\Gamma :=
\cup_{T\in\mathscr{T}} \partial T$ to denote the skeleton of the partition
$\mathscr{T}$, and we define $\Gamma_\text{int} := \Gamma\setminus\partial D$
and $\Gamma_\text{ext} := \Gamma\cap\partial D$. 
\begin{definition}
Let $T^+$ and $T^-$ be two neighboring elements, and let $e\coloneqq\partial T^+\cap
\partial T^-\in\Gamma_\text{int}$. Let ${n}^+$ and ${n}^-$
denote the outward pointing unit normal vector fields on $e$ with respect to $T^+$
and $T^-$. Let $\psi\in H^1(\mathscr{T})$ be a scalar function, and let $\psi^+$ and
$\psi^-$ denote the traces onto $e$ of $\psi\vert_{T^+}$ and $\psi\vert_{T^-}$,
respectively. The average $\avge{\psi}$ and the jump $\jumpe{\psi}$ of $\psi$ on
$e$ are defined by
\begin{equation}
\avge{\psi} \coloneqq \frac{1}{2}(\psi^++\psi^-)
\quad \text{and}\quad
\jumpe{\psi}\coloneqq\psi^+n^++\psi^-n^-\,,
\end{equation}
respectively. Similarly, let $\Psi\in [H^1(\mathscr{T})]^d$ be a vector valued function,
and let $\Psi^+$ and $\Psi^-$ denote the traces onto $e$ of $\Psi\vert_{T^+}$
and $\Psi\vert_{T^-}$, respectively. The average $\avge{\Psi}$ and the jump
$\jumpe{\Psi}$ of $\Psi$ on $e$ are defined by,
\begin{equation}
\avge{\Psi} \coloneqq \frac{1}{2}(\Psi^++\Psi^-)
\quad\text{and}\quad
\jumpe{\Psi} \coloneqq \Psi^+\cdot n^++\Psi^-\cdot n^-\,,
\end{equation}
respectively. Finally, let $\tilde T\in\mathscr{T}$ be such that $\partial
\tilde T\cap \partial D\neq \emptyset$. Let $\tilde e \coloneqq \partial \tilde T\cap
\partial D\subset\Gamma_{\text{ext}}$, and let $\tilde n$ denote the outward
pointing unit vector field on $\tilde e$. Let $\tilde \psi$ and $\tilde \Psi$
denote the traces onto $e$ of $\psi\vert_{\tilde T}$ and $\Psi\vert_{\tilde T}$,
respectively. Using this notation, we extend the definition of jump and average
on $\Gamma_{\text{ext}}$ as
follows:
\begin{equation}
\avge{\psi} \coloneqq \tilde\psi\,,\quad
\jumpe{\psi} \coloneqq \tilde\psi\tilde n\,,\quad
\avge{\Psi} \coloneqq \tilde\Psi\,,\quad
\jumpe{\Psi} \coloneqq \tilde\Psi\cdot\tilde n\,.
\end{equation}
\end{definition}

\subsection{Discontinuous Galerkin discretization of the shape gradient
equation} \label{sec:dgelliptic} Let $\Ch = H^1_0(D,\bbR^d)$ be endowed with the
standard $H^1$-norm, let $E_i$, $i=1,\dots, d$ denote the canonical basis of
$\bbR^d$, and let $g_i\coloneqq \nabla J(\Omega)\cdot E_i$, $i=1,\dots, d$
denote the components of $\nabla J(\Omega)$.  Then, the shape gradient equation
\cref{eq:shapegradienteqn} implies that each $g_i\in H^1_0(D)$ is characterized
by the variational problem
\begin{equation}\label{eq:shapegradienteqnH10}
\int_{D} \big(\nabla g_i \cdot\nabla w + g_i w \big) \,dx = dJ(\Omega; wE_i)
\quad \text{for all }w\in H^1_0(D)\,.
\end{equation}
To derive the dG discretization of \cref{eq:shapegradienteqnH10}, it is
important to highlight that the shape derivative $dJ(\Omega;V)$ generally
depends both on the direction $V$ and its Jacobian $DV$. For this reason, we introduce
the notation
\begin{equation*}
dJ_i(\Omega; w,\nabla w)\coloneqq dJ(\Omega; wE_i)\,.
\end{equation*}
We can now derive the dG discretization of \cref{eq:shapegradienteqnH10}. 
First, we note that each $g_i\in H^1_0(D)$ is the unique minimizer of the
functionals
\begin{equation*}
F_i:H^1_0(D)\ni v\mapsto
0.5\int_{D} \big(\vert \nabla v \vert^2 + v^2 \big)\,dx - dJ_i(\Omega; v, \nabla v)\,.
\end{equation*}
To facilitate consistency of differentiation the presence of element-wise discontinuous
functions, we consider the, so-called, lifting
operator $R: H^1(\mathscr{T})\to [S^p_{\mathscr{T}}]^d$, defined by 
\begin{equation}\label{eq:liftingoperator}
\int_D R(v)\cdot w\dx = -\int_{\Gamma} \jump{v}\cdot \avg{w}\ds
\quad \text{for all }w\in [S^p_{\mathscr{T}}]^d\,,
\end{equation}
we refer to \cite{ArBrCoMa01} for a discussion and properties. 
We define the following discretized version of $F_i$, reading
$F_i^h:S^p_{\mathscr{T}}\to\bbR$,
\begin{equation}\label{eq:discreteshapegradientfunctional}
F_i^h(v)\coloneqq
\frac{1}{2}\int_{D}  \big(\vert\nabla_h v + R(v)\vert^2 + v^2\big) \,dx
+\frac{1}{2}\int_{\Gamma} \sigma \vert \jump{v}\vert^2\ds
- dJ_i(\Omega; v, \nabla_h v + R(v))\,,
\end{equation}
with $\sigma:\Gamma\to\mathbb{R}$ a non-negative function, to be defined
precisely in \cref{sec:penalty}, typically referred to as the
\emph{discontinuity-penalization parameter}.  The minimizer $g_{i,h}\in
S^p_{\mathscr{T}}$ to the discrete functional
\cref{eq:discreteshapegradientfunctional} is characterized by
$dF_i^h(g_{i,h},w)=0$ for all $w\in S^p_{\mathscr{T}}$, and
can thus be computed by solving the finite dimensional variational problem
\begin{equation}\label{LDG-like}
\begin{aligned}
B(g_{i,h},w):=&\ \int_{D}\big( \nabla_h g_{i,h} \cdot \nabla_h w +g_{i,h}w \big)\,\dx\\
&+\int_{\Gamma} \Big(\sigma\jump{g_{i,h}}\cdot\jump{w}
-\jump{w}\cdot \avg{\nabla g_{i,h}}-\jump{g_{i,h}}\cdot \avg{\nabla_h w+R(w)}\Big)\ds\\
=&\ dJ_i(\Omega; w, \nabla_h w + R(w))\qquad\text{for all }w\in S^p_{\mathscr{T}}.
\end{aligned}
\end{equation}

To bypass the requirement of computing the lifting of every test function in
\cref{LDG-like}, it is possible to derive a ``pure'' interior penalty method by
considering the alternative discretized version of $F_i$, reading
\begin{equation}\label{eq:discreteshapegradientfunctional_two}
	\begin{aligned}
\tilde{F}_i^h:S^p_{\mathscr{T}}\ni v\mapsto
\frac{1}{2}&\int_{D}  \big(\vert\nabla_h v \vert^2 + v^2\big) \,\dx
-\int_{\Gamma}\avg{\nabla_h v}\cdot\jump{v}\ds
+\frac{1}{2}\int_{\Gamma} \sigma \vert \jump{v}\vert^2\ds\\
&
- dJ_i(\Omega; v, \nabla_h v + R(v)).
\end{aligned}
\end{equation}
Then, the minimizer $\tilde{g}_{i,h}\in S^p_{\mathscr{T}}$ to the discrete functional
\cref{eq:discreteshapegradientfunctional_two} is is characterized by
$d\tilde{F}_i^h(\tilde{g}_{i,h},w)=0$ for all $w\in S^p_{\mathscr{T}}$, that is, 
\begin{equation}\label{IPDG}
\begin{aligned}
\tilde{B}(g_{i,h},w):=&\ 
\int_{D}\big( \nabla_h \tilde{g}_{i,h} \cdot \nabla_h w +g_{i,h}w \big)\,dx \\
&+\int_{\Gamma} \Big(\sigma\jump{\tilde{g}_{i,h}}\cdot\jump{w}
-\jump{w}\cdot \avg{\nabla \tilde{g}_{i,h}}	
-\jump{\tilde{g}_{i,h}}\cdot \avg{\nabla_h w}\Big)\ds\\
=&\  dJ_i(\Omega; w, \nabla_h w + R(w))\qquad\text{for all }w\in S^p_{\mathscr{T}}.
\end{aligned}
\end{equation}

\begin{example}\label{ex:dJuncostrained}
In \cref{ex:levelsetobjective}, we saw that if $J(\Omega)\coloneqq \int_\Omega f\dx$ for a given scalar function $f\in W^{1,1}(D)$, then 
\begin{equation}\label{eq:shapederivativelevelsetfunctional}
dJ(\Omega; V) = \int_\Omega \big( \nabla f\cdot V + f\Div{V}\big)\dx
\quad \text{for }V\in W^{1,\infty}(D;\bbR^d)\,.
\end{equation}
The extension of 
formula \cref{eq:shapederivativelevelsetfunctional} to vector fields of the
form $V=wE_i$ with $w\in S^p_{\mathscr{T}}$ reads
\begin{align}\nonumber
dJ_i(\Omega; w, \nabla w+ R(w)) = 
&\int_\Omega (\nabla f\cdot E_i) w\dx
+ \int_{\Omega} f\nabla_h w\cdot E_i\dx+ \int_{\Omega} \Pi_\Omega fR(w)\cdot E_i\dx\\
\label{eq:dJuncostrained}
=&\int_\Omega (\nabla f\cdot E_i) w\dx
+ \int_{\Omega} f\nabla_h w\cdot E_i \dx
- \int_{\Gamma\cap \bar{\Omega}}\! \jump{w}\avg{\Pi_\Omega fE_i}\ds\,.
\end{align}
\end{example}

As we shall see in the numerical experiments in \cref{sec:Numericalexperiments},
general polytopic elements will be constructed via \emph{agglomeration} from a
finer simplicial subdivision of the hold-all domain $D$. This is particularly
pertinent in the present context of shape optimization, as we shall see below,
due to the flexibility it offers in capturing small geometrical features with
relatively coarse elements.  At the same time, however, constructing polytopic
meshes via agglomeration may result into very irregular element shapes. To
ensure the stability of the dG method, therefore, we resort into sophisticated
choice of the discontinuity-penalization parameter $\sigma$, as derived and
justified in \cite{CaDoGeHo17,CaDoGe21}.

\subsection{The discontinuity-penalization parameter}\label{sec:penalty}
For completeness, we now give a precise formula for the
discontinuity-penalization function $\sigma$ appearing in the interior penalty
discontinuous Galerkin (IP-dG) formulations above.  The stability and error
analysis of the IP-dG method under this choice of penalization, as well as a
detailed discussion on the practical relevance of this choice, can be found in
\cite{CaDoGeHo17,CaDoGe21}.

The construction of a suitable discontinuity-penalization function $\sigma$
begins with the concept of (finite simplicial) covering. Each polytopic element
$T\in \mathscr{T}$ can be covered by a finite family
$\mathcal{K}_{T}:=\{K_j\}_{j=1}^{m}$ of $m\in\bbN$ simplices $K_j$,
i.e., we have $T\subset \cup_{K_j\in \mathcal{K}_{T}} K_j$.
We refer to $\mathcal{K}_{T}$ as a \emph{covering} of $T$.  Of course, many
coverings with different cardinalities exist. For instance, any simplicial
partition  of $T$ is a covering; equally, coverings with overlapping simplices
also exist. We shall denote by $\mathbb{K}_{T}$ the set of all such coverings.

For a generic domain $\omega\subset\mathbb{R}^s$, $s\in \{1,\dots, d\}$, let
$h_\omega$, $\rho_\omega$  and $|\omega|$ denote the diameter, the inscribed
radius, and the $sD$-volume of $\omega$, respectively. 

\begin{definition}
We say that an element
$T$ is \emph{$p$-coverable} if there exists at least one covering
$\mathcal{K}_{T}$ of $T$ such that
\[
\max_{x\in T}\{\min_{z\in \partial K_j}|{x}-{z}|\}\le \rho_{K_j}/(8p)^2,
\]
with $|\cdot|$ denoting the Euclidean distance, for some constants $c_{sh}>0$ and $c'_{sh}>0$, such that
  $c_{sh}^{-1}  h_{T}\le h_{K_j}\le c_{sh}  h_{T}$, and
$(c'_{sh})^{-1}  \rho_{T}\le \rho_{K_j}\le c'_{sh}  \rho_{T}$,  $c'_{sh}>0$,
for all $K_j\in \Ck_T$.
\end{definition}
In other words, an element is
$p$-coverable if there exists a covering $\mathcal{K}_T$ comprising simplices
each with similar shape-regularity to the original polytopic element $T$, which
cover $T$ within a distance at most $\rho_{K_j}/(8p)^2$ each, away from the
element's boundary $\partial T$. In \cite{CaDoGe21}, a considerably weaker concept
of $p$-coverability is used allowing, in particular, $K_j$ to be general
\emph{curved} prisms. 

Let $e\subset\Gamma_{\rm {int}}$ denote a common face shared by two neighboring
polytopic elements $T^+,T^-\in\mathscr{T}$. Let $\underline{\mathbb{K}}_T^e$
denote the set of all simplices that are contained in $T$ and that contain $e$
as part of a face. We define the restriction of the discontinuity-penalization
parameter $\sigma$ on $e$ by
\begin{equation}\label{eq:sigma_skeleton}
\sigma|_e:=C_{\sigma} \max_{T\in\{T^+,T^-\}}\bigg\{ \min \Big\{
\frac{|T|}{\sup_{K\in\underline{\mathbb{K}}_T^e} |K|},C_{\rm cov}(T)
\Big\} \frac{ p_T^2|e|}{|T|} \bigg\}\,
\end{equation}
for a computable constant $C_{\sigma}>0$, depending on $c_{sh}$ and on $c'_{sh}$,
and where we set $C_{\rm cov}(T):=p_T^{2(d-1)}$ if
$T$ is $p$-coverable, or $C_{\rm cov}(T):=\infty$ if $T$ is not $p$-coverable.
The definition of the discontinuity-penalization parameter $\sigma$ can be
naturally extended to boundary faces $e\subset\partial D\cap T$ for a
$T\in\mathscr{T}$ as follows
\begin{equation}\label{eq:sigma_bdry}
\sigma|_e:=C_{\sigma} \min \Big\{
\frac{|T|}{\sup_{K\in\underline{\mathbb{K}}_T^e} |K|},C_{\rm cov}(T)
\Big\} \frac{ p_T^2|e|}{|T|} \,.
\end{equation}

\begin{remark}
It is possible to extend this definition for extreme local element size and
polynomial degree variations in the context of Robust IP-dG methods
\cite{DoGe22}; we refrain from doing so in the present work, preferring to focus
on the shape-optimization aspects.
\end{remark}

The definition of $\sigma$ given in formulas
\cref{eq:sigma_skeleton,eq:sigma_bdry} is designed for maximum generality.
However, it may appear at first sight to be complicated to be implemented in
practice. If we additionally assume that each element $T\in\mathscr{T}$ is
star-shaped with respect to a ball $B({x}_0,ch_T)$ for some $c<1$ and
$x_0\in T$, then it is sufficient to select 
\[
\sigma|_e:=C_{\sigma} \max_{T\in\{T_1,T_2\}}\big\{ p_T^2|e| h_T^{-d}\big\}
\quad\text{and}\quad
\sigma|_e:=C_{\sigma} p_T^2|e| h_T^{-d}\,,
\]
respectively.

\subsection{Recovery}\label{sec:recovery} The requirement $V\in
[W^{1,\infty}(D)]^d$ in \cref{eq:levelseteqn} is sufficient for the
well-posedness of the respective PDE problem. However, \eqref{LDG-like} and
\eqref{IPDG} produce element-wise \emph{discontinuous} approximations
$g_{i,h}\in  S^p_{\mathscr{T}}$.  A natural and minimally invasive approach to
address this is to perform local recovery of the vector field
$(g_{i,h})_{i=1}^d$ into $H^1_0(D)$. For this purpose, it is sufficient to
perform a simple nodal averaging approach upon constructing a suitable
sub-division of the polytopic elements into $d$-dimensional simplices. 

More specifically, we consider $\tilde{\mathscr{T}}$ to be a refinement of
$\mathscr{T}$ into simplices, meaning that each polytopic element
$T\in\mathscr{T}$ is subdivided into $r_T$ simplices $\tau_m\subset T$ with
$m=1,\dots,r_T$, so that $\bar{T}=\cup_{m=1}^{r_T}\bar{\tau}_m$. (The same, or
different simplicial subdivision of each polytopic element can be used for
quadrature computations during assembly, see also \cite{CaDoGeHo17,DoGeKa21}.) Related
to the simplicial mesh $\tilde{\mathscr{T}}$, we consider the discontinuous
element-wise polynomial space \begin{equation} S^p_{\tilde{\mathscr{T}}} := \{v
\in L^2(D)\,:v|_\tau \in {P}_p(\tau) \text{ for all }
\tau\in\bar{\mathscr{T}}\}\,, \end{equation} for which we have, trivially, $
S^p_{\mathscr{T}} \subset S^p_{\tilde{\mathscr{T}}} $. We also denote by
$\tilde{\Gamma}$ the skeleton of $\tilde{\mathscr{T}}$.
Further, on $S^p_{\tilde{\mathscr{T}}}$, for each Lagrange node $x_j$ of each
element $\tau$, globally numbered via $j=1,\dots, {\rm
dim}S^p_{\tilde{\mathscr{T}}} $, we  consider the nodal element patch $ \omega
(x_j):=\{ \tau \in \tilde{\mathscr{T}}: x_j\in\bar{\tau}\}$, with cardinality
denoted by ${\rm card}\big(\omega(x_j)\big)$. We define the nodal averaging
recovery operator $\mathcal{E}:S^p_{\tilde{\mathscr{T}}} \to H^1(D)\cap
S^p_{\tilde{\mathscr{T}}}$ by
\begin{equation*}
\mathcal{E}(v)(x_j) \coloneqq \frac{1}{{\rm
card}\big(\omega(x_j)\big)}\sum_{\tau\in\omega(x_j)}v|_{\tau}(x_j)
\end{equation*}
for $v\in S^p_{\tilde{\mathscr{T}}} $ and $x_j\in D$, while $\mathcal{E}(v)(x_j)
=0$ for all $x_j\in\partial D$. Note that if a Lagrangian node $x_j$ lies in the
interior of a simplex $\tau$, then $\mathcal{E}(w)(x_j) =w(x_j)$, because in
this case $\omega(x_j)=\tau$.

With this, by solving \eqref{LDG-like} or \eqref{IPDG}, injecting the
numerical solution into $S^p_{\tilde{\mathscr{T}}}$, and applying recovery
operator $\Ce$, we obtain an approximate shape gradient $\nabla J(\Omega)_h\in
[S^p_{\tilde{\mathscr{T}}} \cap H^1_0(D)]^d$. Note that, since $\nabla
J(\Omega)_{\tilde h}$ is continuous and element-wise polynomial, we also have
$\nabla J(\Omega)_{\tilde h}\in [W^{1,\infty}(D)]^d$.  

\begin{remark}
A completely analogous construction is possible when the fine mesh
$\tilde{\mathscr{T}}$ comprises box-type elements, i.e., quadrilaterals when
$d=2$ or hexahedra for $d=3$. The only modification is that the respective
recovery space is then the space of element-wise mapped tensor-product elements,
as is standard in conforming finite elements. 
\end{remark}

\subsection{Runge-Kutta discontinuous Galerkin discretization of the level-set
equation} \label{sec:dghyperbolic}

Having constructed a recovered $\nabla J(\Omega)_h\in [W^{1,\infty}(D)]^d$, we
consider the level-set equation \cref{eq:levelseteqn} with $V=-\nabla
J(\Omega)_h$. To approximate its solution $\varphi$ on $D$, we employ a
Runge-Kutta discontinuous Galerkin (RKdG) method for first order equations on a
simplicial (or box-type) mesh $\mathscr{T}^{fine}$, i.e., a refinement into
simplicial (or box-type)
elements of the polytopic mesh $\mathscr{T}$. For instance, we may take
$\mathscr{T}^{fine}=\tilde{\mathscr{T}}$ as in \cref{sec:recovery}, or we could even
consider $\mathscr{T}^{fine}$ to be a refinement of $\tilde{\mathscr{T}}$.
For notational simplicity, we shall, henceforth, take
$\mathscr{T}^{fine}=\tilde{\mathscr{T}}$, stressing that this is not necessary
for the overall framework to work.

We note that, employing a spatial discretization of  \cref{eq:levelseteqn} on
$\tilde{\mathscr{T}}$ (or a refinement thereof!), is \emph{not} necessarily
detrimental to the overall complexity of the method. This is due to the explicit
time-stepping nature of RKdG, since no global linear system is required to be
solved on each time-step, and since quadratures take place on simplicial
subdivisions. {Indeed, the main computational cost
(associated with the time-stepping component of an explicit RKdG method)
consists in assembling elemental mass matrices and
solving the resulting block-diagonal linear systems. 
Using $L^2$-orthogonal bases, which are available for simplices, 
the mass matrices become diagonal, making this computational costs negligible.
Of course, it is by all means possible to employ general polytopic meshes for the
discretization of the level-set equation \cite{CaDoGeHo16}. 
In this case, 
assembly of elemental mass matrices typically involves quadrature
computations on a simplicial subdivision of the polytopic mesh, in which case the
computational cost is the same as computing directly on the simplicial subdivision.
For large problems, this can be significantly accelerated via GPU implementations \cite{DoGeKa21}.
Using polytopic meshes may lead to
fewer but larger blocks compared to a discretization on the underlying simplicial mesh.
However, the resulting polytopic elemental mass matrices are generally not diagonal, and so the fast
solution aspect may be compromised.
We also note that the spatial accuracy is generally driven by elements' diameters, and so solutions computed on
polytopic meshes may be less accurate than counterparts obtained on their simplicial subdivision.
Therefore, given the comparable computational cost and potentially higher accuracy, we opt for solving the level-set equation on the simplicial subdivision $\tilde{\mathscr{T}}$.}

For a generic element $\tau\in\tilde{\mathscr{T}}$, we denote by $\partial
\tau^-$ the inflow region of $T$, that is
\begin{equation}
	\partial \tau^-\coloneqq\{ x\in\partial \tau :V\cdot n < 0\}\,.
\end{equation}
Note that $\partial \tau^-\subset\Gamma_\text{int}$ because $V=-\nabla
J(\Omega)_{\tilde h}$ vanishes on $\partial D$ due to our choice of $\Ch$ in
\cref{sec:dgelliptic}.

For a sufficiently smooth function $w:D\to\mathbb{R}$, we define the so-called
\textit{upwind jump} at $x\in\partial \tau^-\subset\Gamma_\text{int}$ to be
\begin{equation}
	\upwind{w}(x) \coloneqq w^+(x) - w^-(x)\,,
	\quad \text{where} \quad
	w(x)^\pm \coloneqq \lim_{\epsilon\to 0^+}w(x\pm\epsilon V)\,.
\end{equation}

Let assume that the initial condition $\phi$ lies in $S^p_{\tilde{\mathscr{T}}}$
(otherwise, $\phi$ can be just replaced by its $L^2$-projection 
$\Pi_\Omega(\phi) \in S^p_{\tilde{\mathscr{T}}}$). The spatially discrete
discontinuous Galerkin method for the level-set problem \eqref{eq:levelseteqn}
on $\tilde{\mathscr{T}}$ reads: find $\varphi_{\tilde{h}}\in
S^p_{\tilde{\mathscr{T}}}$ such that $\varphi_{\tilde{h}}(0)=\phi$ in $D$ and
\begin{equation}\label{eq:levelsetDG}
	\int_D (\varphi_{\tilde h})_t w \,dx + \sum_{\tau\in\tilde{\mathscr{T}}}\int_\tau (V\cdot\nabla\varphi_{\tilde{h}}) w \, dx 
	- \sum_{\tau\in\tilde{\mathscr{T}}}\int_{\partial \tau^- }(V\cdot n)\upwind{\varphi_h}w^+\,ds = 
	0 \quad\text{for all }w\in 	S^p_{\tilde{\mathscr{T}}} \,.
\end{equation}
The spatially discrete dG method \eqref{eq:levelsetDG} does not include upwind
boundary terms because $V\vert_{\partial D}=0$. 

Finally, we discretize  \eqref{eq:levelsetDG} with respect to the $t$-variable
using explicit Runge-Kutta methods that have favorable properties when combined
with spatial dG discretizations. We briefly review these methods and their
properties here for the sake of completeness.

Selecting a local elemental basis $\{\psi_\tau^i\}_{i=1}^{ \dim(P_p(\tau))}$, we
express the elemental mass and ``stiffness'' matrices as
$M_\tau:=(m^\tau_{ij})_{i,j=1}^{\dim(P_p(\tau))}$ and
$K_\tau:=(k^\tau_{ij})_{i,j=1}^{\dim(P_p(\tau))}$, where
\[
m^\tau_{ij}:=\int_\tau \psi_T^j \psi_T^i \,dx
\quad\text{and}\quad k^\tau_{ij}:=\int_\tau (V\cdot\nabla \psi_T^j)  \psi_T^i \, dx 
- \int_{\partial \tau^- \setminus \partial D^-}(V\cdot n)\upwind{\psi_T^j}(\psi_T^i)^+\,ds,
\]
respectively. We note that, upon selecting orthogonal elemental basis, e.g., the
so-called Dubiner system \cite{KaSh05}, $M_\tau$ becomes diagonal, and thereby
trivially invertible.

Let $\varphi_{\tilde h}\vert_\tau(t) :=
\sum_{j=1}^{\dim(P_p(\tau))}\Phi_\tau^j(t)\psi_\tau^j$, let
$\Phi_\tau(t):=(\Phi_\tau^1(t),\dots, \Phi_\tau^{ \dim(P_p(\tau))}(t))$, and let
$L_\tau:=-M_\tau^{-1}K_\tau$. With this notation, we can write
\eqref{eq:levelsetDG} on each $\tau\in\tilde{\mathscr{T}}$ as
\begin{equation}\label{semi_DG}
\frac{d}{dt}\Phi_\tau(t)=L_\tau \Phi_\tau(t),\qquad t\in(0,\epsilon),
\end{equation}
On \eqref{semi_DG}, we apply a standard stable Runge-Kutta method, such as
such as Heun's method, or TVB-type methods \cite{CoShu98}, selecting carefully the appropriate CFL restrictions.

\subsection{Polytopic dG discretization of the level-set shape optimization
method} \label{sec:dgpolytopic}

We can now describe how to exploit the polytopic dG methods ability to simply
and accurately resolve domain boundaries (the zero-level set) in the computation
of the shape derivative and of eventual PDE constraints and adjoint equations.
We first recall from \cref{sec:dghyperbolic} that the level-set function
$\varphi$ is approximated on the simplicial mesh $\tilde{\mathscr{T}}$ using a
Runge-Kutta dG method. Since the approximation $\varphi_{\tilde h}$ is
discontinous, we postprocess it with the recovery strategy described in
\cref{sec:recovery} to ensure that the implicitly defined domain
$\Omega_{\Ce(\varphi_{\tilde h})}$ is well defined.  Then, we identify the
simpleces that intersect the zero level set by tracking the change of sign of
$\Ce(\varphi_{\tilde h})$ on the mesh nodes. For example, in 2D we mark the
$i^{\text{th}}$ triangle if
\begin{equation*}
\left|\frac{\Ce(\varphi_{\tilde h})(x_{i,1})}{|\Ce(\varphi_{\tilde h})(x_{i,1})|}
+\frac{\Ce(\varphi_{\tilde h})(x_{i,2})}{|\Ce(\varphi_{\tilde h})(x_{i,2})|}
+\frac{\Ce(\varphi_{\tilde h})(x_{i,3})}{|\Ce(\varphi_{\tilde h})(x_{i,3})|}\right|<3\,,
\end{equation*}
where $x_{i,1}, x_{i,2}$, and $x_{i,3}$ denote the triangle vertices. Here, we
tacitly assume that if $\Ce(\varphi_{\tilde h})$ is a piecewise polynomial of
degree higher than one, then $\tilde{\mathscr{T}}$ is sufficiently fine so that,
e.g. in 2D, the recovered level-set function $\Ce(\varphi_{\tilde h})$ does not
change sign twice along any mesh edge.

\begin{figure}[htb!]
\centering
\includegraphics[width=0.7\linewidth]{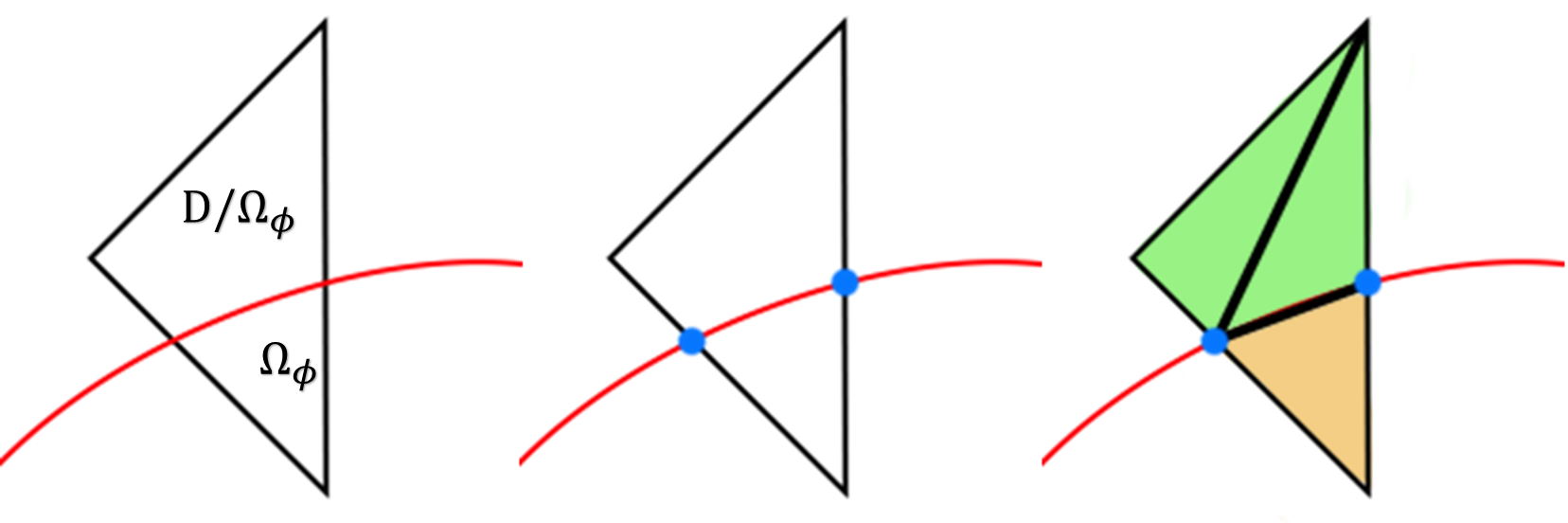}
\caption{Example of local refinement method in 2D. Upon each edge that is
intersected by the zero-level set $\{\phi = 0\}$ (red curved line), we
interpolate the level set function and determine its zeros, where we append
new nodes (black dots) used to split the triangle in (at most) three
subtriangles. {To obtain a more accurate representation of
the zero-level set, one could first subdivide elements that intersect the zero-level
set and then apply the local refinement method to newly generated elements that
intersect the zero-level set, possibly in combination with a higher-order moving-mesh
shape optimization method \cite{PaWeFa18}.}}
\label{fig:localrefinement}
\end{figure}

After marking the simplices, we locally refine them to better resolve the
zero-level set, thus creating a refined partition $\mathscr{T}^\text{ref}\supset
\tilde{\mathscr{T}}$ (see \cref{fig:localrefinement} for an example in 2D).
Finally, we agglomerate the simplices in $\mathscr{T}^\text{ref}$ to form the
polytopic subdivision $\mathscr{T}^\text{agg} \subset \mathscr{T}^\text{ref}$,
which can be used to compute the approximation of the shape gradient as
discussed in \cref{sec:dgelliptic}.  For example, in the numerical experiments
presented in \cref{sec:Numericalexperiments}, we create $\mathscr{T}^\text{agg}$
using the following simple graph partitioning method based on $k$-means
clustering \cite{MJ67}. First, we determine the barycenters of the triangles in
$\mathscr{T}^\text{ref}$ and interpret them as the set of observations
$\{x_j\in\mathbb{R}^{2}\,: j=1,\dots, n_\text{ref}\}$, where $n_\text{ref}$
denotes the number of simplices in $\mathscr{T}^\text{ref}$. Then, we split the
set of observations in two depending on the sign of $\phi(x_j)$ and introduce
$k^+\in\mathbb{N}$ and $k^-\in\mathbb{N}$ distributed sets $C_i^+$ and $C_i^-$.
We assign each barycenter $x_j$ with $\phi(x_j)>0$ to a distributed set $C_i^+$
by minimizing the quantity
\begin{equation*}
\sum_{i=1}^{k^+}\sum_{x\in C^+_i}
\left\Vert{x-\frac{1}{|C^+_i|}\sum_{y\in C^+_i}y}\right\Vert^2\,,
\end{equation*}
where $\vert C_i^+\vert$ denotes the cardinality of $C_i^+$ and $\Vert
\cdot\Vert$ denotes the Euclidean distance. Barycenters $x_j$ with $\phi(x_j)<0$
are assigned to $C_i^-$ in a similar fashion.  Finally, we agglomerate triangles
with barycenters in the same distributed set to form the $(k^++k^-)$ polytopic
elements in $\mathscr{T}^\text{agg}$. Examples of
$\tilde{\mathscr{T}},\mathscr{T}^\text{ref}\;\text{and}\;\mathscr{T}^\text{agg}$
can be seen in both \cref{fig:unconstrained_flow} and
\cref{fig:constrained_flow}.

We conclude this section by briefly explaining how to define and integrate
polynomial basis functions on arbitrarily shaped polygonal elements.  For a more
detailed treatment of the polytopic dG finite element method, we refer to \cite{CaDoGeHo14,CaDoGeHo17}.

The construction of basis functions on polytopic elements is based upon defining
tensor-product polynomial spaces on each element using tightly fitting bounding
boxes. Each bounding box is chosen such that the most extreme vertices of the
polytope $T$ align with the edges of the bounding box so that we acquire a
minimum bounding box $B_T$. Upon this bounding box $B_T$, we can construct a set
of basis functions by restricting onto $T$ the tensor-product Legendre
polynomials spanning $P_p(T)$.

To assemble the Galerkin matrices, we take advantage of the underlying
simplicial mesh $\mathscr{T}^\text{ref}$ used to form the polytopes $T \in
\mathscr{T}^\text{agg}$. More specifically, we define quadrature schemes on
polytopes by simply combining standard quadrature rules for simpleces. 

\section{Numerical experiments}\label{sec:Numericalexperiments}

We now demonstrate the good numerical performance of the level-set shape
optimization method based on polytopic discontinuous Galerkin discretizations.

\subsection{Optimization algorithm}\label{sec:numexp_algo} As explained in
\cref{sec:shapeopttheory}, starting from an initial guess $\Omega_{\phi^0}$, we
can generate a sequence of functions $\phi^k$ such that $J(\Omega_{\phi^k})$ is
decreasing in a steepest-descent with line-search fashion \cite{NoWr06}. More
specifically, we set $\phi^{k+1} \coloneqq \varphi^k(t^k)$ for $k\geq 0$, where:
$\varphi^k$ is a numerical solution to
\begin{equation}\label{eq:alg_levelseteqn}
\partial_t \varphi^k - \nabla J(\Omega_{\phi^k})\cdot \nabla\varphi^k = 0
\quad \text{for }t\in(0,\epsilon_k)\,, \quad \varphi^k(0) = \phi^{k}\,,
\end{equation}
$t^k\in (0,\epsilon_k)$ is a suitable stopping time, and $\nabla
J(\Omega_{\phi^k})$ is a numerical solution to
\begin{equation}\label{eq:alg_shapegradienteqn}
(\nabla J(\Omega_{\phi^k}), V)_\Ch = dJ(\Omega_{\phi^k}; V)\quad \text{for all } V\in \Ch\,,
\end{equation}
for a chosen Hilbert space $\Ch$. In the following numerical experiments, we use
$\Ch = H^1_0(D;\bbR^d)$, but other choices are possible. The final time
$\epsilon_k$ must be chosen such that $\varphi^k(t)$ exists for every $t\in (0,
\epsilon_k)$. Note that both $\epsilon_k$ and the time step $\Delta t_k$ used to
approximate $\varphi^k$ depend on the $W^{1,\infty}(D;\bbR^d)$-norm of $\nabla
J(\Omega_{\phi^k})$, which possibly changes at each iteration $k$. In the
following numerical experiments, we select $\Delta t_k$ to satisfy the
CFL-condition and perform at most $M\in\bbN$ time steps, thus implicitly setting
$\epsilon_k=M\Delta t_k$. The stopping time $t^k$ is selected in a backtracking
fashion as the largest in the set $\{{5}m\Delta t_k\,: m\in \bbN, 0<{5}m\leq M\}$
such that the Armijo condition
\begin{equation}\label{eq:alg_Armijo}
J(\Omega_{\varphi^k(m\Delta t_k)})\leq
J(\Omega_{\varphi^k(0)}) - c (m\Delta t_k) \Vert\nabla J(\Omega_{\phi^k})\Vert^2_\Ch
\end{equation}
is satisfied, where the constant $c\in (0,1)$ is set to $c=0.01$. The
optimization algorithm is stopped after 80 iterations, or if if the Armijo
condition \cref{eq:alg_Armijo} is not satisfied after one time step, that is, if
\begin{equation*}
J(\Omega_{\varphi^k(\Delta t_k)}) > 
J(\Omega_{\varphi^k(0)}) - c \Delta t_k \Vert\nabla J(\Omega_{\phi^k})\Vert^2_\Ch\,.
\end{equation*}

\begin{figure}[htb!]
\centering
\includegraphics[width=\linewidth]{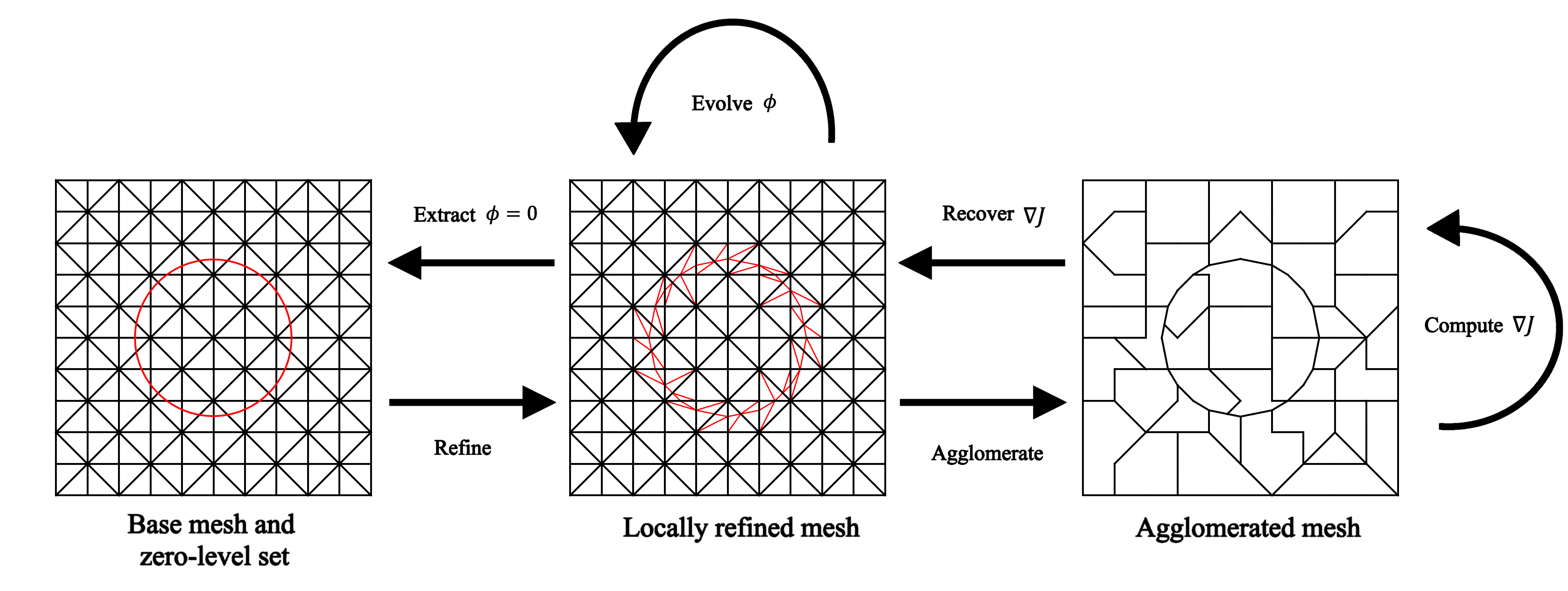}
\caption{Workflow of the agglomerated
polytopic dG level-set shape optimization method. A base triangular mesh {(left)} is locally
refined to better fit the zero-level set $\phi=0$. The resulting mesh {(middle)}
is subsequently agglomerated {into a polytopic mesh (right)}
to compute an approximate solution to the shape gradient equation
\cref{eq:alg_shapegradienteqn}. The solution $\nabla J$ is then
{recovered (via interpolation and nodal averaging)} on the
locally refined triangular mesh,
where the level-set equation \cref{eq:alg_levelseteqn}
is solved using a time-stepping scheme. When the time-stepping scheme stops,
the {base mesh is }locally refined mesh to better fit the zero-level set of the updated
level-set function $\phi$. This procedure is repeated until the convergence or stopping
criteria are met.}
\label{fig:unconstrained_flow}
\end{figure}

\begin{figure}[htb!]
\centering
\includegraphics[width=\linewidth]{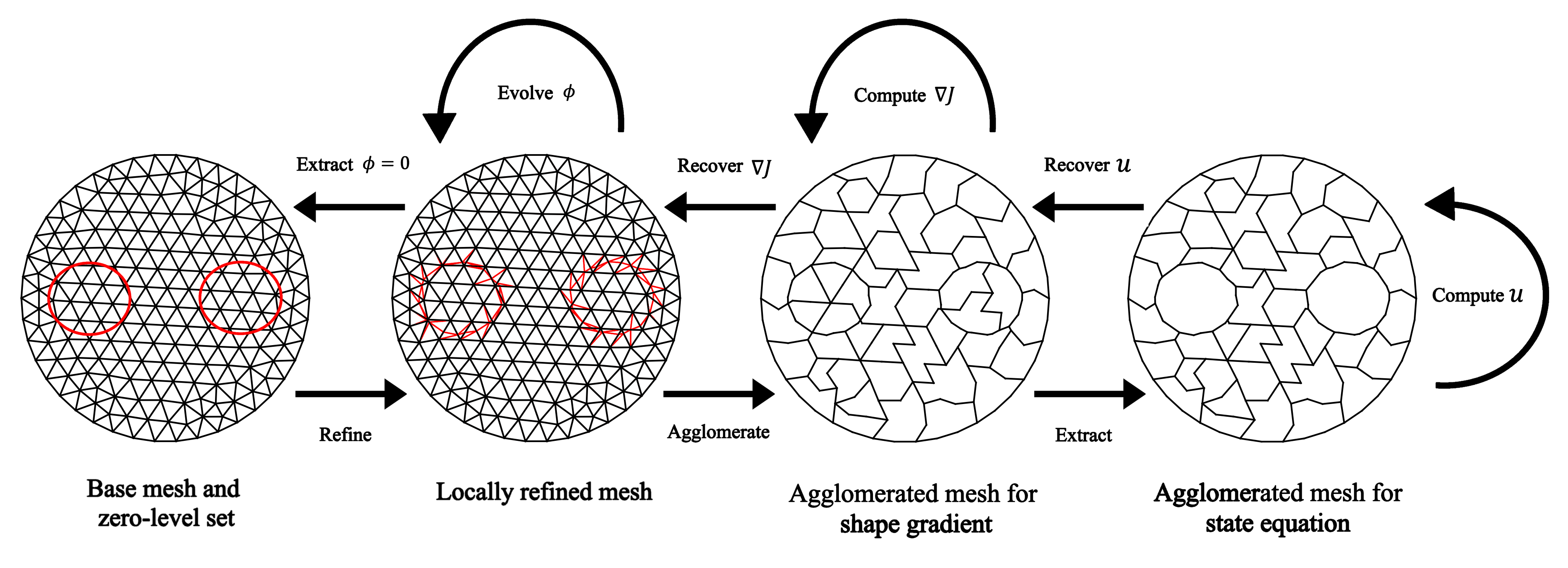}
\caption{Extension of the agglomerated polytopic dG level-set shape optimization
workflow presented in \cref{fig:unconstrained_flow} to PDE-constrained shape
optimization problems. {The workflow begins as in \cref{fig:unconstrained_flow}.}
To compute the shape gradient $\nabla J$, we first
extract the state equation's computational domain and compute the state solution
{$u$. The state solution is then recovered (via interpolation of $u$, $\nabla u$, and the lift $R(u)$)
on the polytopic mesh used to compute the shape gradient $\nabla J$. Then, the workflow continues
as in \cref{fig:unconstrained_flow}.}}
\label{fig:constrained_flow}
\end{figure}

The workflow of the resulting optimization algorithm for uncostrained problems
is sketched in \cref{fig:unconstrained_flow}. Its extension to PDE-costrained
problems is depicted in \cref{fig:constrained_flow}. We anticipate that in these
numerical experiments, as well as in every other we performed (not
reported here), the level-set shape optimization algorithm based on the
discontinuous Galerkin method did not require a reinitialition of the level set
function $\phi$, a procedure that is commonly invoked in the literature
\cite{AlDaJo21,AlDaFr14, Ad20}. This is in agreement with the findings in \cite{MaReCh06}
in the case the level-set equation in divergence form for divergence-free velocity vectors.

\subsection{Unconstrained test case} \label{sec:unconstrainedtestcase}

We consider the unconstrained shape functional $J(\Omega)\coloneqq\int_\Omega
f\dx$, where the bivariate function $f$ is chosen as 
\begin{equation*}
f(x)\coloneqq f(x_1,x_2)=\sqrt[4]{((x-0.7)^2 + y^2)((x+0.7)^2 + y^2)}-0.6\,.
\end{equation*}

To solve the level-set equation \cref{eq:alg_levelseteqn}, we partition the
hold-all domain $D$ into 1800 triangles, discretize in space using RKdG
with spatial basis functions of degree $p = 2$ and timestepping via Heun's
method with time step $\Delta t_k=1/2600$, with at most $M={150}$
time steps per iteration. Choosing a very small time step
$\Delta t_k$ is not a requirement of the dG method employed. This choice is
simply due to the desire to generate ``smooth'' convergence histories for
plotting purposes and to provide a ``smallest step'' stopping criteria for the
backtracking algorithm. To solve the shape gradient equation
\cref{eq:alg_shapegradienteqn}, we agglomerate the locally refined triangular mesh
(cf.~\cref{fig:unconstrained_flow}) into 200 polytopic elements and employ
polynomials of degree $p=2$. {\Cref{tab:shapegradientconvergence}
reports a numerical validation of the polytopic discontinuous Galerkin approximation
of the shape gradient}. The computation of the shape derivate $dJ$ in the
discontinuous Galerkin framework is discussed in \cref{ex:dJuncostrained},
see  \cref{eq:dJuncostrained}. To assemble the various Galerkin matrices, we
employ quadrature rules of order four both on triangles and on edges, which compute
the involved integrals exactly (up to round-off errors).

\begin{table}[htbp]
\footnotesize
\caption{
Numerical validation of the polytopic discontinous Galerkin approximation of the
shape gradient $\nabla J$ in the unconstrained test case. The computational domain is a square of
edge length 2 and the zero-level set is a circle of radius 0.52 (both centered at the
origin, as depicted in white in the most left image in \cref{fig:shapeevolutionunconstrained}).
In the table, the index $\mathrm{i}$ denotes
the level of refinement, $\mathrm{N}$ denotes the number of polytopic elements, and
$\mathrm{err}_l$ and $\mathrm{err}_q$ denote the $L^2$-error for affine and quadratic polygonal finite
element approximations, respectively. The approximate reference solution has been computed
using standard continuous quadratic finite elements on a simplicial mesh with 4944 elements.
We observe that the approximate algebraic converge rates
computed using the formula
$\mathrm{rate}(i)=\log(\mathrm{err}(\mathrm{i})/\mathrm{err}(\mathrm{i}-1))/\log(\mathrm{N}(\mathrm{i})/\mathrm{N}(\mathrm{i}-1))$
are close to $-1$ and $-1.5$, which are the expected rates for finite element
approximations based on affine and quadratic polynomials, respectively.}
\label{tab:shapegradientconvergence}
\begin{center}
  \begin{tabular}{|c|r|c|c|c|c|} \hline
  $\mathrm{i}$ & $\mathrm{N}$ &  $\mathrm{err}_l$& \bf $\mathrm{rate}_l$ &  $\mathrm{err}_q$& \bf $\mathrm{rate}_q$ \\ \hline
    1 & 36 & 1.34e-3 & - & 3.09e-4 & -\\
    2 & 144 & 4.06e-4 & -0.86 & 6.70e-5& -1.10\\ 
    3 & 576 & 1.35e-4 & -0.79&7.25e-6 &-1.60\\ 
    4 & 2304 & 3.86e-5& -0.91 &1.04e-6 &-1.40\\ \hline
  \end{tabular}
\end{center}
\end{table}

\begin{figure}[htb!]
\centering
\includegraphics[width=0.24\linewidth]{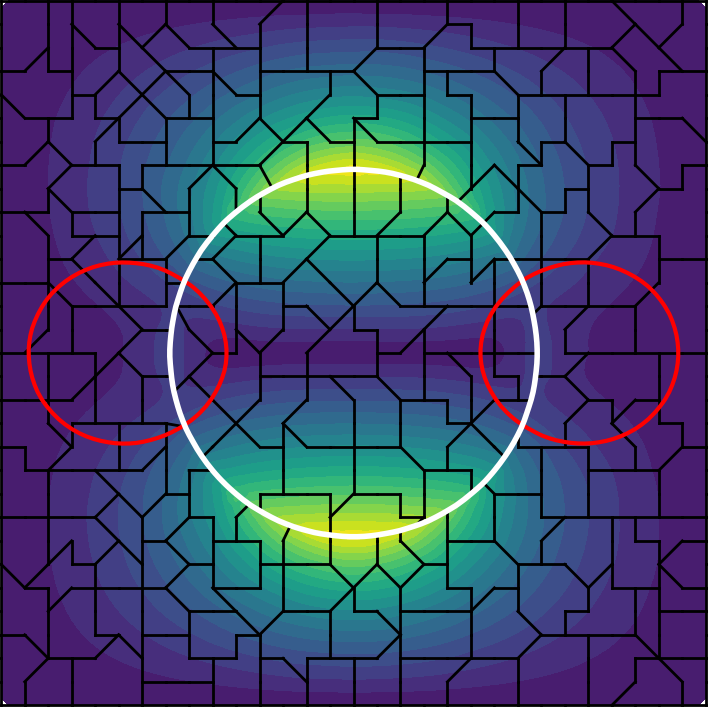}
\includegraphics[width=0.24\linewidth]{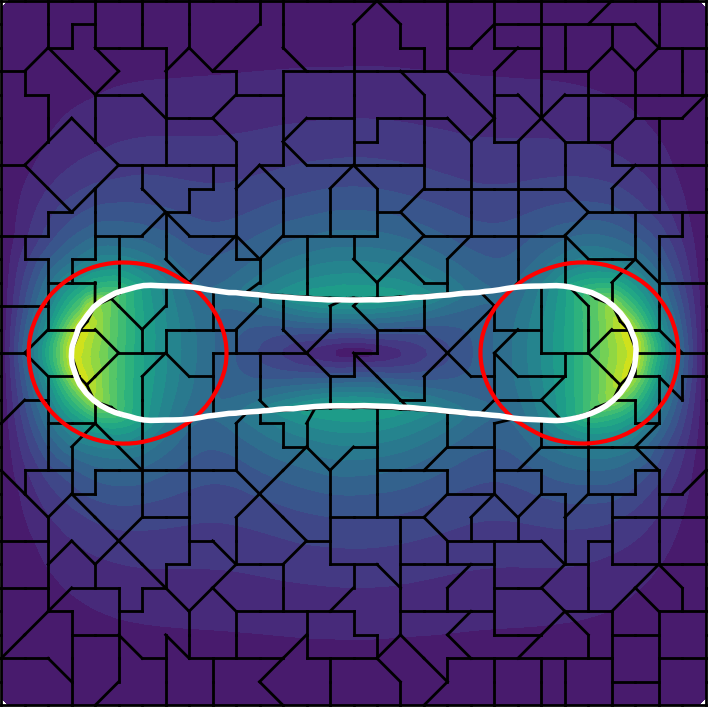}
\includegraphics[width=0.24\linewidth]{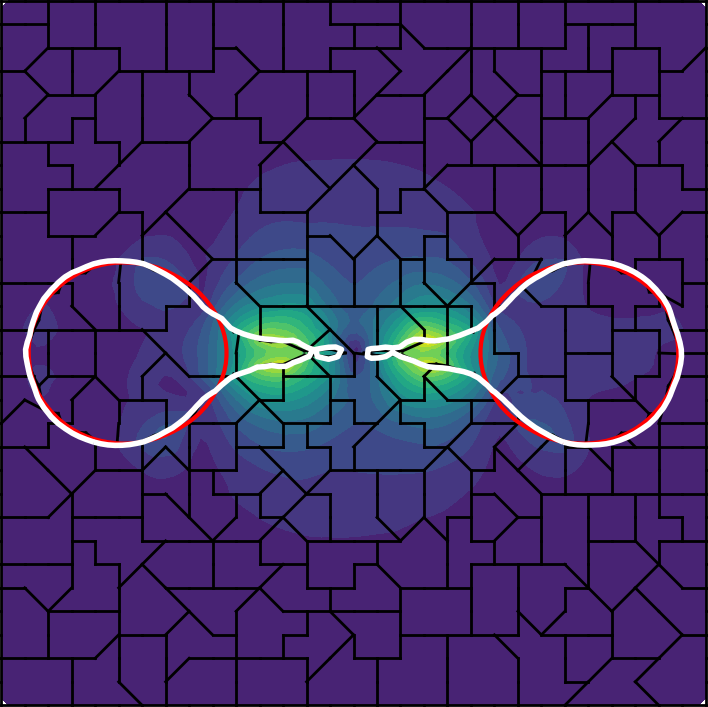}
\includegraphics[width=0.24\linewidth]{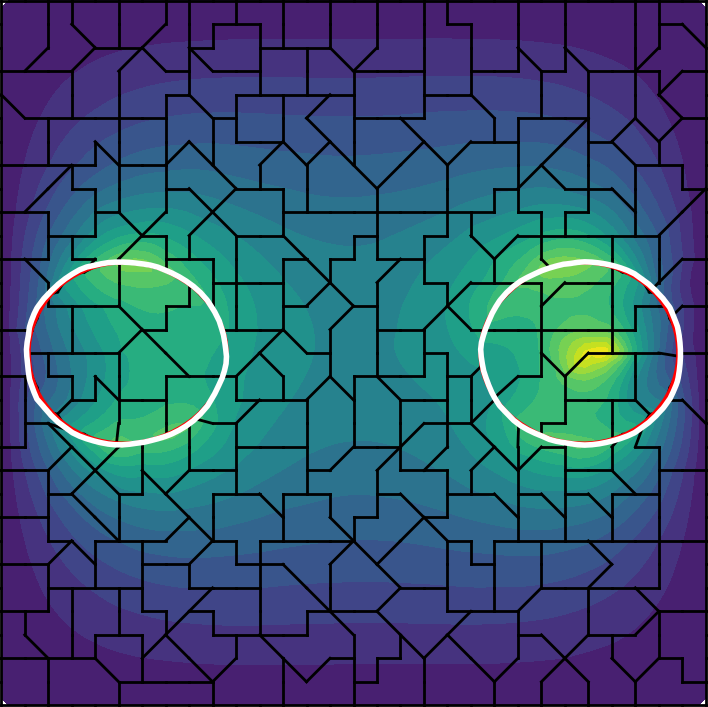}\\
\begin{tikzpicture}[overlay]
\node at (-4.8,0) {$n =0$, $t = 0$};
\node at (-1.6,0) {$n = 9$, $t = 0.52$};
\node at (1.6,0) {$n = 14$, $t = 0.81$};
\node at (5.05,0) {$n = 41$, $t = 1.45$};
\end{tikzpicture}
\caption{Evolution of shapes for the unconstrained test case.
{The variables $n$ and $t$
denote the number of gradient evaluations and the pseudo-time, respectively.}
The initial domain,
a disc {(drawn in white)}, elongates and eventually splits into two domains before converging
to the optimal configuration {(drawn in red)}. The contour plots represent the
normalized contours of the Euclidean norm of the gradient $\nabla J$.
Note that the polytopic mesh changes at each iteration.}
\label{fig:shapeevolutionunconstrained}
\end{figure}

\begin{figure}[htb!]
\centering
\includegraphics[width=0.6\linewidth]{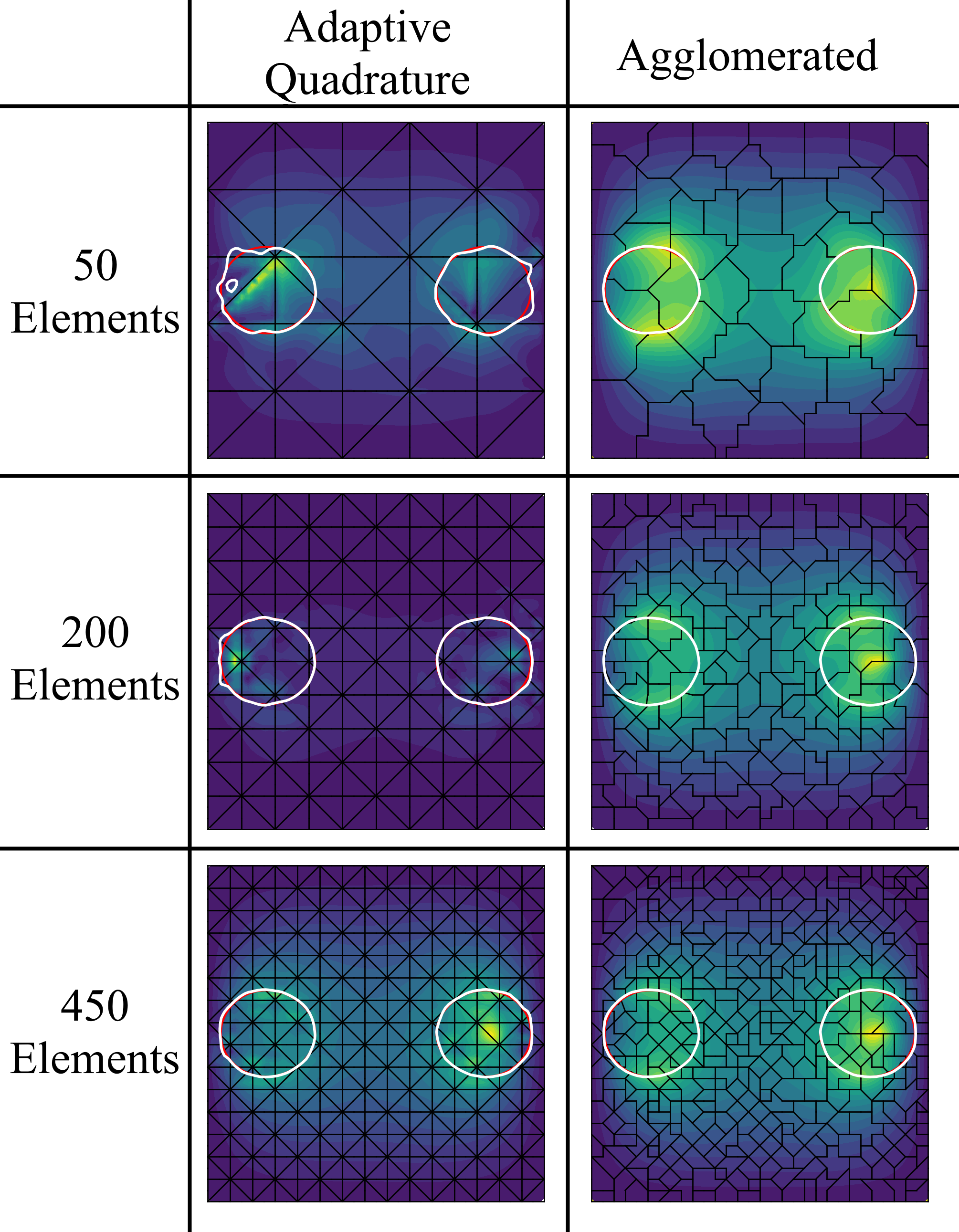}
\caption{{Qualitative comparison
(based on the unconstrained test case)}
between the agglomerated polytopic dG approach and a more
classical approach based on triangular elements with adaptive quadrature.
{The optimal solution is the union of two disjoint disk-shaped domains depicted in red.
The optimized shapes are depicted in white. The background shows the meshes employed and} the
normalized contours of the Euclidean norm of the gradient $\nabla J$.
{We observe that, as expected, both approaches retrieve increasingly accurate solutions
as the the number of mesh elements increases. The agglomerated polytopic dG approach
with 50 elements is qualitatively comparable, if not better, than
the adaptive quadrature approach with 450 elements.}}
\label{fig:qualitative_comparison_unconstrained}
\end{figure}

\begin{figure}[htb!]
\centering
\includegraphics[width=0.49\linewidth, trim=20 40 20 20, clip]{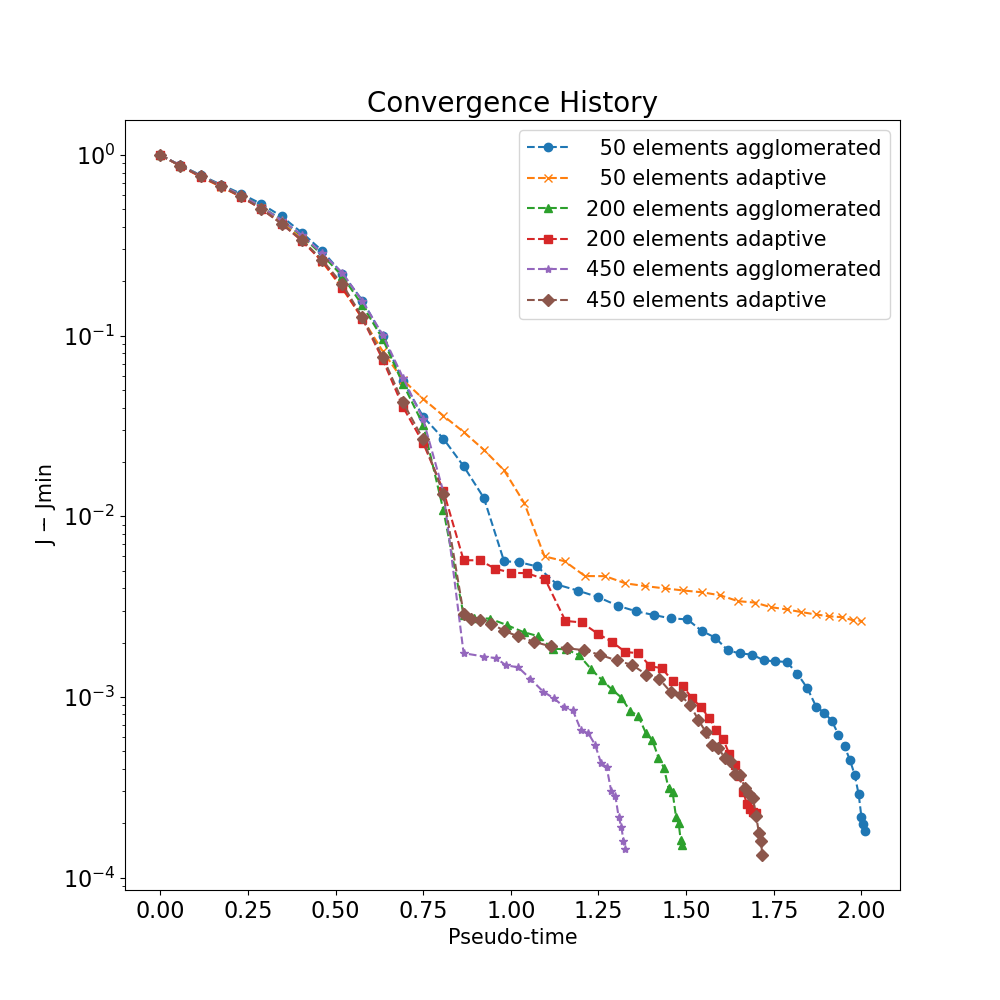}
\includegraphics[width=0.49\linewidth, trim=20 40 20 20, clip]{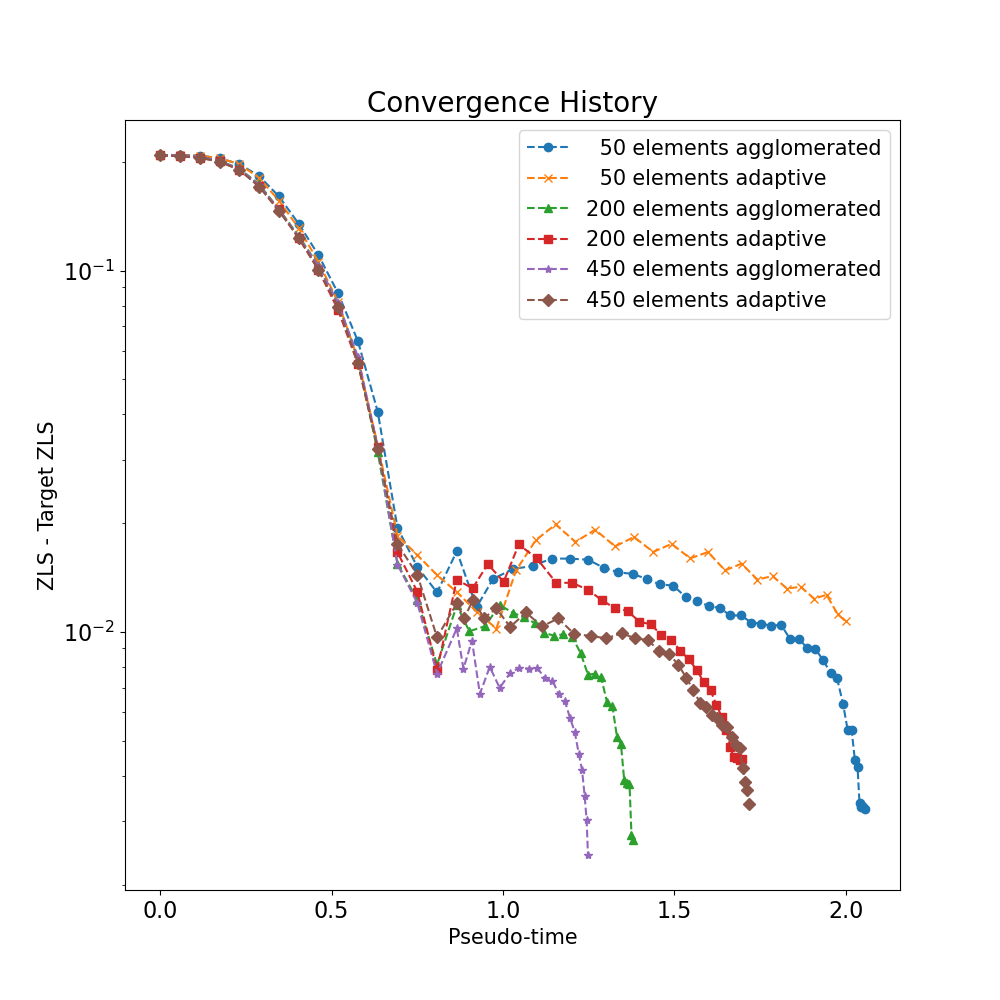}
\caption{Quantitative comparison (based on the unconstrained test case)
between the agglomerated polytopic dG approach and a more
classical approach based on triangular elements with adaptive quadrature.
The figure on the left depicts the pseudo-time evolution of the difference between the
objective value $J$ as shapes are optimized and (an accurate approximation of) the
objective value at the optimal solution.
The figure on the right depicts the pseudo-time evolution of an approximate distance between
the zero-level set of the optimized shapes and the zero-level set of the optimal shape.
Markers denote gradient evaluations.
We observe that the agglomerated polytopic dG generally converges more rapidly and accurately
than the adaptive quadrature approach.}
\label{fig:quantitative_comparison_unconstrained}
\end{figure}

Some shape iterates, starting from a disc of radius 0.51 centered at the origin,
are displayed in \cref{fig:shapeevolutionunconstrained}.
{The normalized pseudo-time evolution of the objective function is
depicted in \cref{fig:quantitative_comparison_unconstrained}}
We observe that the
algorithm retrieves a good approximation of the optimal (and not connected)
shape. To assess the advantage of the agglomerated polytopic dG approach to
compute shape derivatives, we compare it with a more classical dG approach based
on triangular meshes. For the latter, we also employ an adaptive quadrature
based on locally refined meshes to evaluate the shape derivative formula $dJ$
accurately on triangles that intersect the level-set boundary $\phi = 0$. For
both methods, we solve the level-set equation \cref{eq:alg_levelseteqn} on a fine
underlying simplicial mesh as in
the previous numerical experiment.
\Cref{fig:qualitative_comparison_unconstrained,fig:quantitative_comparison_unconstrained}
summarizes the qualitative and quantitative comparison between the two
approaches when using {50, 200, and 450} triangular or polytopic elements.
{The quantitative comparison is performed both in terms of the objective functions $J$ and
a numerical approximation of the distance between zero-level sets defined as
$\max_{x\in \{\varphi=0\}} \min_{y\in\{f=0\}} \vert x - y\vert$.}
We observe
that the adaptive nature of the agglomerated polytopic approach allows it to
achieve accurate results also with coarser meshes.

\subsection{PDE-constrained test case}\label{sec:constrainedtestcase}
We consider the PDE-constrained objective function $ J(\Omega, u)\coloneqq
\int_\Omega \nabla u\cdot\nabla u+\eta^2\dx$, where $u\in H^1(\Omega)$ is the
weak solution to
\begin{equation}\label{eq:stateconstraint} 
-\Delta u = 0 \quad \text{in }\Omega\,,
\quad u=f\quad\text{on }\partial \Omega\,.
\end{equation}
This PDE-constrained optimization problem stems from the class of Bernoulli free
boundary value problems. In these problems, a part $\Gamma_{\mathrm{free}}$ of
the boundary $\Gamma$ is moved so that the solution $u$ to the state problem
\cref{eq:stateconstraint} satisfies simultaneously both the Dirichlet and
Neumann boundary conditions $u=f$ and $\partial_n u = \eta$ on
$\Gamma_{\mathrm{free}}$ \cite{EpHa12}. Specifically, we set the hold-all $D$ to
be a unit disc centered at the origin. The boundary $\partial D$ of $D$
coincides with the fixed part $\Gamma_{\mathrm{fixed}}\coloneqq
\Gamma\setminus\Gamma_{\mathrm{free}}$ of the boundary $\Gamma$. The free
boundary $\Gamma_{\mathrm{free}}$ is the boundary of subdomains compactly
embedded in $D$.  Furthermore, we set $f=-1$ on $\Gamma_{\mathrm{fixed}}$, $f=0$
on $\Gamma_{\mathrm{free}}$, and $\eta = 1/(-0.55\ln{(0.55)})$, so that the
optimal shape is an annulus with internal radius equal to 0.55. This data has
been generated starting for the fact that, if $\Omega$ is an annulus and $f$ is
constant on each connected component of $\partial \Omega$, then the solution to
\cref{eq:stateconstraint} is of the form $u(x)=a\log(\vert x \vert) + b$, where
$a,b\in\bbR$ are two constants that depend on the $\Omega$ and $f$.

In this setting, the shape derivative of $J$ in a sufficiently regular
direction $V$ is \cite{HiPa15}:
\begin{equation}\label{eq:dJconstrainedstd}
dJ(\Omega, u; V)=\int_\Omega \nabla u\cdot
(\Div(V)-DV^\top - DV)\nabla u + \eta^2 \Div(V)\dx.
\end{equation}
To extend formula \cref{eq:dJconstrainedstd} to directions of the form $V=wE_i$,
where $w\in S^p_{\mathscr{T}}$ and $E_i$ $i=1,\dots, d$ denotes the canonical
basis of $\bbR^d$, we note that $\Div(V)= \nabla w\cdot E_i$, and $DV = E_i
(\nabla w)^\top$. Therefore,
\begin{align*}
dJ_i(\Omega,u;w,\nabla w+R(w)) = &\sum_{T\in\mathscr{T}}\int_{T\cap\Omega} \Big(
(\eta^2+\nabla u\cdot\nabla u) (\nabla w + R(w))\cdot E_i\\
&-\nabla u\cdot ((\nabla w+R(w)) E_i^\top
+ E_i (\nabla w+R(w))^\top)\nabla u\Big)\dx\,.
\end{align*}
Using the definition of the lifting operator \eqref{eq:liftingoperator}, 
since $\eta^2\in S^p_\mathscr{T}$, we can rewrite
\begin{equation*}
\sum_{T\in\mathscr{T}}\int_{T\cap\Omega} 
(\eta^2+\nabla u\cdot\nabla u) R(w)\cdot E_i\dx
=- \int_{\Gamma_{\text{int}}\cap \bar{\Omega}}
\jump{w}\avg{(\eta^2+\nabla u\cdot\nabla u)E_i}\ds\,.
\end{equation*}
Similarly,
\begin{align*}
\sum_{T\in\mathscr{T}}\int_{T\cap\Omega} 
\nabla u\cdot (R(w) E_i^\top + E_i R(w)^\top)\nabla u\dx
&= \sum_{T\in\mathscr{T}}\int_{T\cap\Omega} 
2(\nabla u\cdot E_i)\nabla u\cdot R(w)\dx \\
&=- \int_{\Gamma_{\text{int}}\cap \bar{\Omega}}
\jump{w}\avg{2(\nabla u\cdot E_i)\nabla u}\ds\,.
\end{align*}
Therefore, we conclude that
\begin{align*}
dJ_i(\Omega,u;w,\nabla w+R(w)) = &\sum_{T\in\mathscr{T}}\int_{T\cap\Omega} 
(\eta^2+\nabla u\cdot\nabla u)E_i\cdot \nabla w
-2(\nabla u\cdot E_i)\nabla u\cdot\nabla w \dx\\
&-\int_{\Gamma_{\text{int}}\cap \bar{\Omega}}
\jump{w}\avg{(\eta^2+\nabla u\cdot\nabla u)E_i-2(\nabla u\cdot E_i)\nabla u}\ds\,.
\end{align*}

Similarly as in \cref{sec:unconstrainedtestcase}, to solve the level-set
equation \cref{eq:alg_levelseteqn}, we partition the hold-all domain $D$ into
{2085} triangles, discretize in space using a classical dG method with basis
functions of degree $p=2$, and employ the Heun's method (with time step $\Delta
t_k=1/1000$) for time-stepping. Per iteration, we perform at most $M={100}$ time steps. To solve the shape gradient equation
\cref{eq:alg_shapegradienteqn} and the state equation
\eqref{eq:stateconstraint}, we agglomerate the locally refined triangular mesh
(cf.  \cref{fig:unconstrained_flow}) into 200 polytopic elements and employ
polynomials of degree $p=2$. To assemble the various Galerkin matrices, we
employ a quadrature rule of order 4 on triangles, an of order 4 on edges, which
computes the involved integrals exactly (up to round-off errors). Note that, as
displayed in \cref{fig:unconstrained_flow}, the state constraint is solved by
extracting its computational domain instead of pursuing the so-called ``Ersatz
material'' approach \cite{AlDaJo21}.

\begin{figure}[htb!]
\centering
\centering
\includegraphics[width=0.24\linewidth]{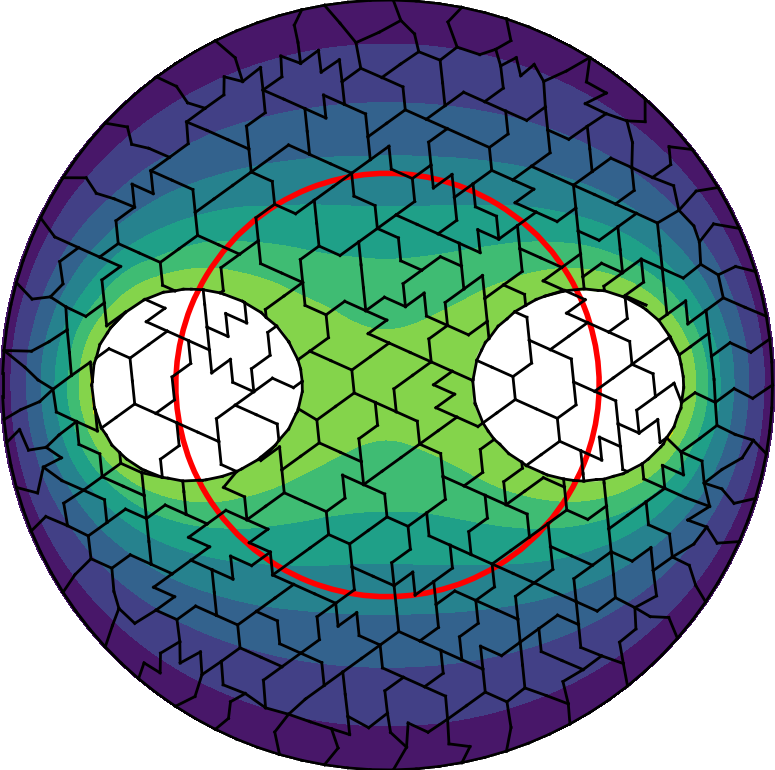}
\includegraphics[width=0.24\linewidth]{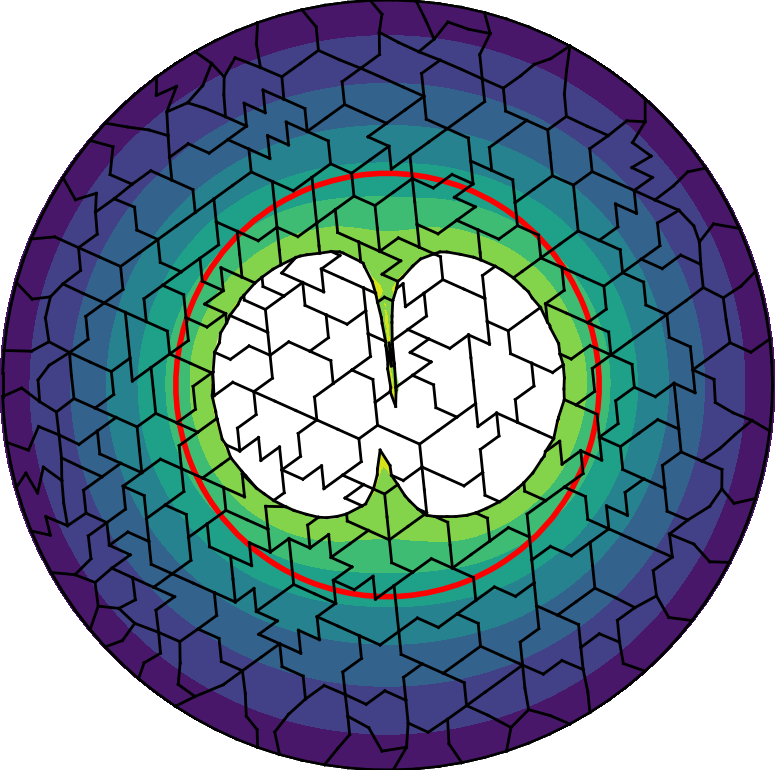}
\includegraphics[width=0.24\linewidth]{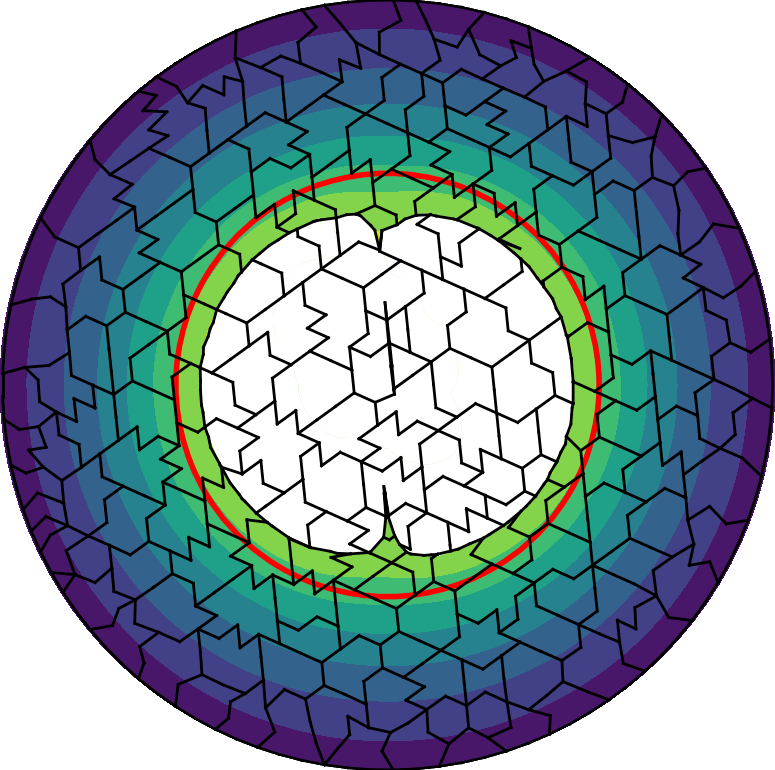}
\includegraphics[width=0.24\linewidth]{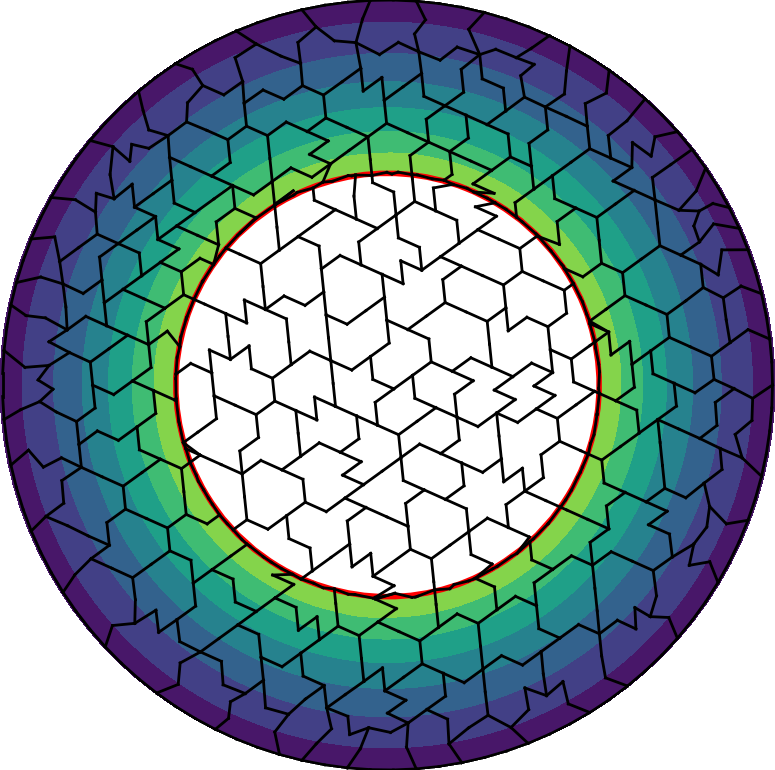}\\
\begin{tikzpicture}[overlay]
\node at (-4.8,0) {$n = 0$, $t = 0$};
\node at (-1.5,0) {$n = 5$, $t = 0.5$};
\node at (1.6,0) {$n = 7$, $t = 0.7 $};
\node at (5.05,0) {$n = 11$, $t = 1.03$};
\end{tikzpicture}
\caption{Evolution of shapes
{for the PDE-constrained test case. The variables $n$ and $t$
denote the number of gradient evaluations and the pseudo-time, respectively.
The two holes in the initial
domain merge into one and then expand,} converging
to the optimal configuration {(depicted in red)}.
The contour plots represent the
normalized contours of the Euclidean norm of the gradient $\nabla J$.
Note that the polytopic mesh changes at each iteration.}
\label{fig:shapeevolutionconstrained}
\end{figure}

\begin{figure}[htb!]
\centering
\includegraphics[width=0.24\linewidth]{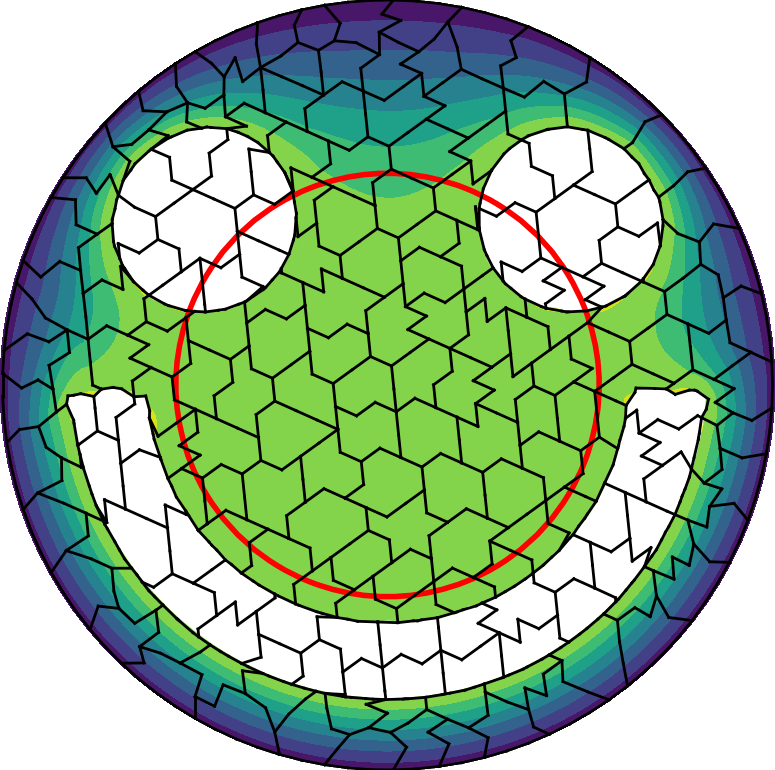}
\includegraphics[width=0.24\linewidth]{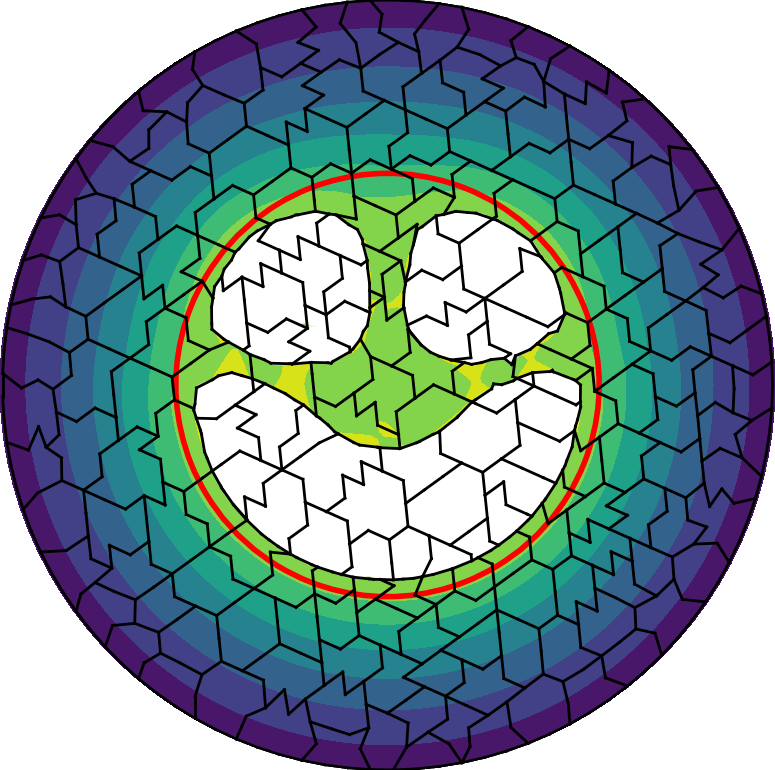}
\includegraphics[width=0.24\linewidth]{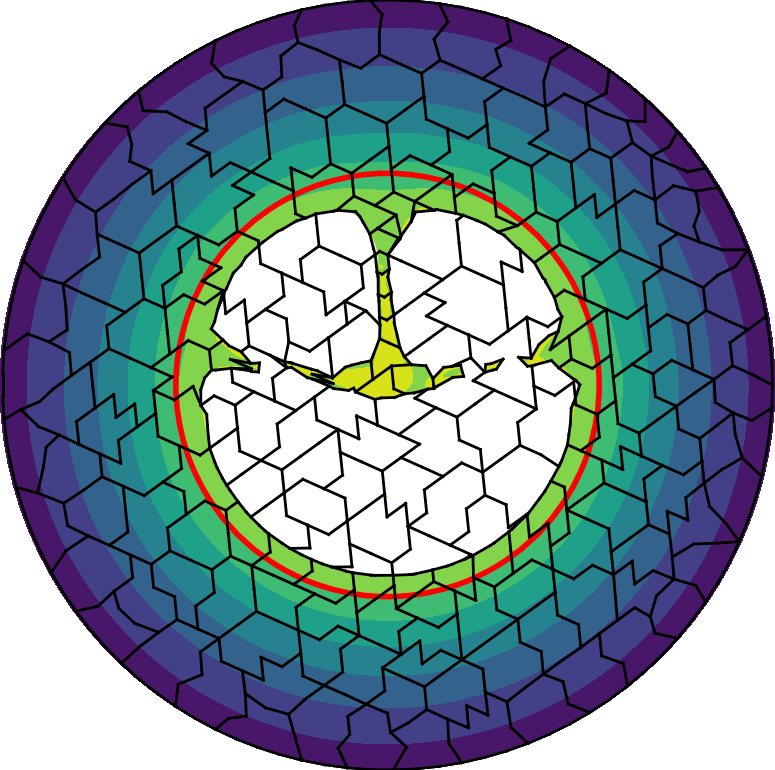}
\includegraphics[width=0.24\linewidth]{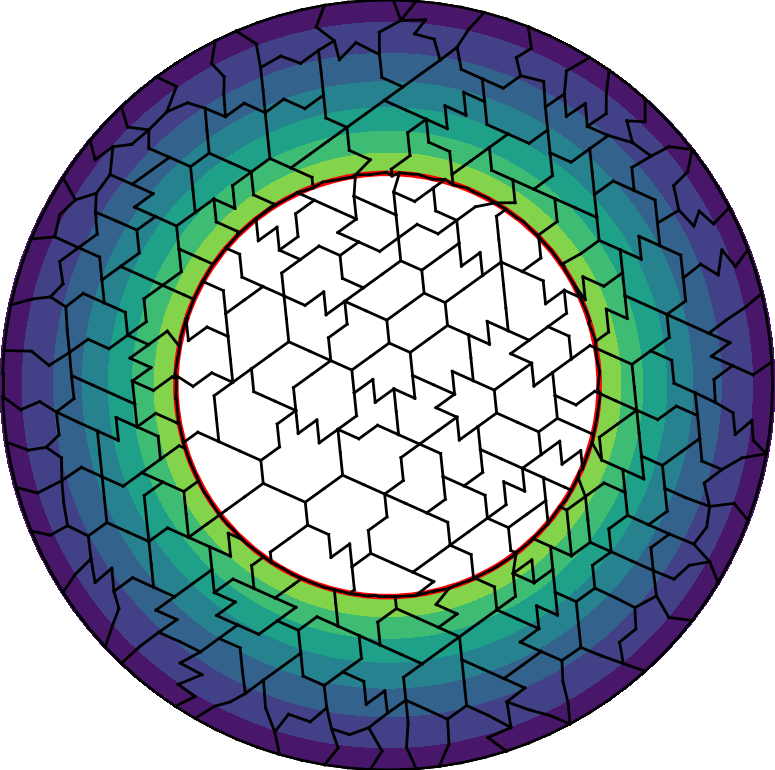}\\
\begin{tikzpicture}[overlay]
\node at (-4.8,0) {$n = 0$, $t = 0$};
\node at (-1.5,0) {$n = 5$, $t = 0.5$};
\node at (1.6,0) {$n = 7$, $t = 0.7$};
\node at (5.05,0) {$n = 13$, $t = 1.1$};
\end{tikzpicture}
\caption{Evolution of shapes {for the PDE-constrained test case with
a different initial guess. The variables $n$ and $t$
denote the number of gradient evaluations and the pseudo-time, respectively.}.
The shape iterates happily converge to the optimal annular domain {(depicted in red)}.}
\label{fig:shapeevolutionconstrained_smiley}
\end{figure}

\begin{figure}[htb!]
\centering
\includegraphics[width=0.49\linewidth]{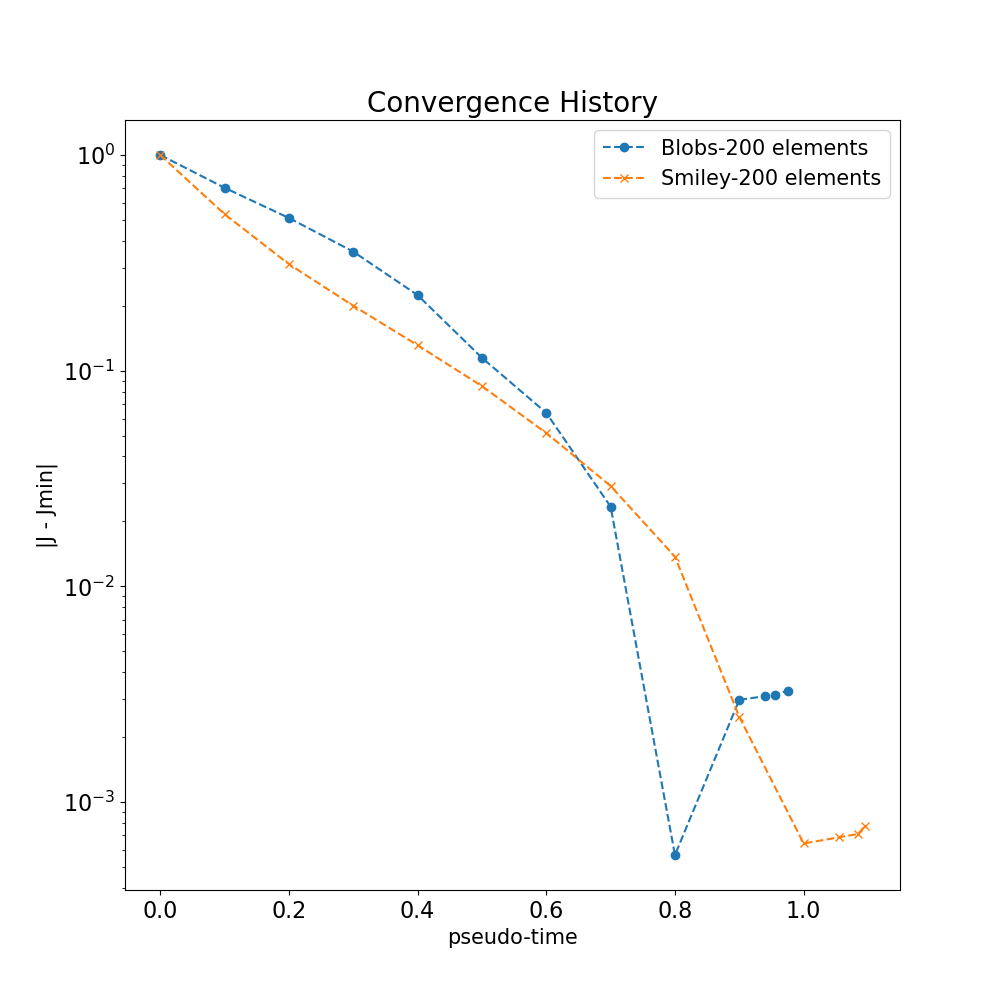}
\includegraphics[width=0.49\linewidth]{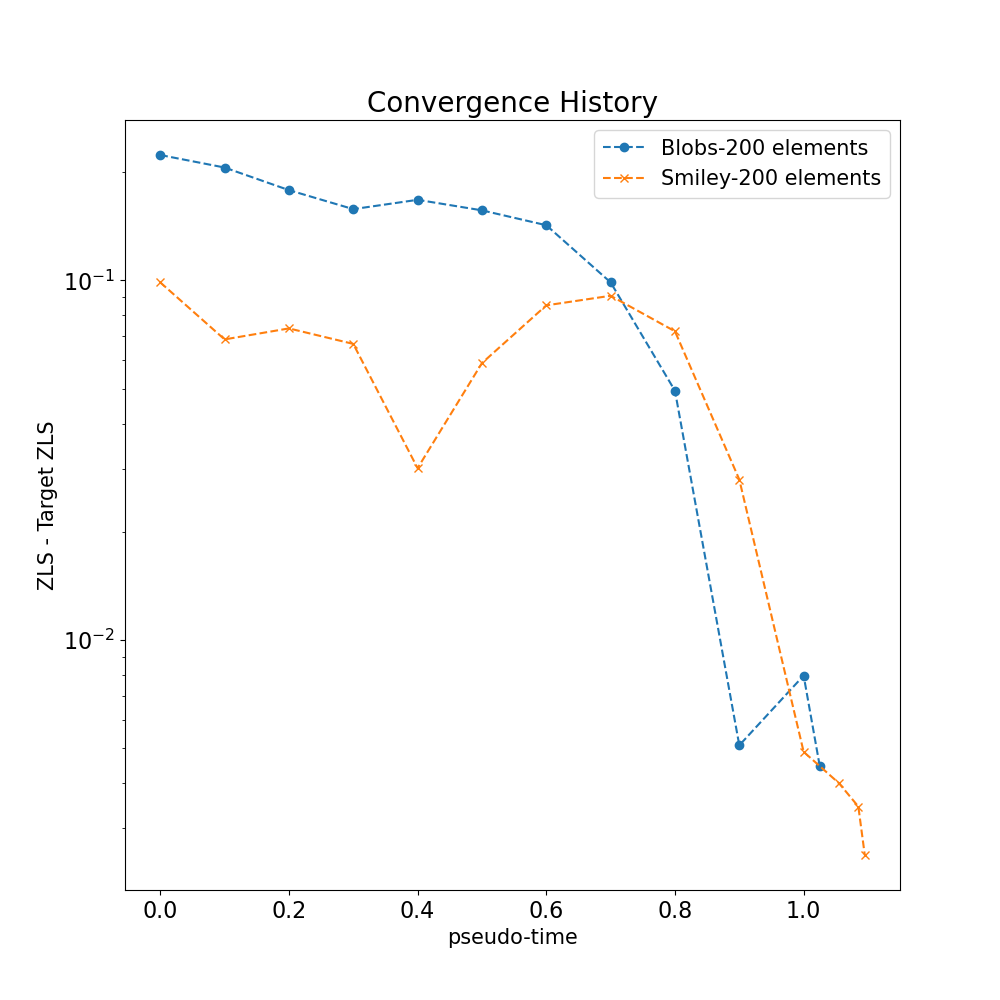}
\caption{Convergence history for the PDE-constrained
test case with different initial guesses: one as in \cref{fig:shapeevolutionconstrained}
and one as in \cref{fig:shapeevolutionconstrained_smiley}.
The figure on the left depicts the pseudo-time evolution of the difference between the
objective value $J$ as shapes are optimized and (an accurate approximation of) the
objective value at the optimal solution.
The figure on the right depicts the pseudo-time evolution of an approximate distance between
the zero-level set of the optimized shapes and the zero-level set of the optimal shape.
Markers denote gradient evaluations.
We observe that both convergence measures are significantly reduced as optimization proceeds.
The convergence lines in the plot on the left plateau due to numerical approximations errors leading
to objective values $J$ fictitiously smaller than objective value at the optimal solution.}
\label{fig:shapeevolutionconstrained_convergencehistory}
\end{figure}

Some shape iterates starting from
\begin{equation*}
\phi(x)=\phi(x_1, x_2)=\sqrt[4]{((x-0.6)^2 + y^2)((x+0.6)^2 + y^2)}-{0.55}\,,
\end{equation*}
are displayed in \cref{fig:shapeevolutionconstrained}. 
We observe that the
algorithm retrieves a good approximation of the optimal shape. To demonstrate
the generality of the approach, \cref{fig:shapeevolutionconstrained_smiley}
displays some shape iterates starting from a different initial guess.
{The convergence histories for both initial guesses are
displayed in \cref{fig:shapeevolutionconstrained_convergencehistory}.}

\section{Conclusions}\label{sec:conclusions}

We have presented a new level-set shape optimization method based on polytopic
discontinuous Galerkin discretizations.
More specifically, the new
shape optimization method exploits the geometric flexibility and robustness of
polytopic dG methods, which allows computing with higher order polynomials
defined on arbitrarily shaped polytopes that resolve the zero-level set
accurately with no need for mesh adjustments. The method also takes advantage of
the efficiency and stability of explicit Runge-Kutta discontinuous Galerkin
methods to solve the level-set equation on a fine simplicial meshes, thus
bypassing the need to reinitialize the level-set function.

We provided a detailed description of the method derivation and explained how to
modify the evaluation of shape derivatives to compute consistent shape gradient 
approximations using discontinuous Galerkin methods. The numerical experiments
evidenced the good performance of the method on unconstrained and
PDE-constrained test cases. An open-source Python implementation of the
method is available on the Zenodo archive at \cite{FeGePa24}.

Thanks to its flexibility in terms of element shapes, polynomial order,
treatment of hanging nodes, stability, the method presented here is perfectly
positioned for the development of \emph{a posteriori} $hp$-adaptive shape optimization
methods. {The above developments can be combined with
GPU-accelerated implementations of $hp$-version polytopic IPdG methods
\cite{DoGeKa21} to enable extremely fast solution of the RKdG step.}

{We note that, although the numerical experiments presented are based on two-dimensional shape optimization examples, three-dimensional problems can be treated using the same methodology, without the need of any further considerations. The favorable computational complexity observed is expected to be retained in the context of three-dimensional setting also. Finally, given the availability of interior penalty dG methods for various classes of PDE problems, the treatment of for more complex PDE problems is by all means possible.}

\bibliographystyle{siamplain}
\bibliography{refs}
\end{document}